\documentclass{amsproc}
\usepackage{amssymb}
\newtheorem{lemma}{Lemma}[section]
\newtheorem{corollary}[lemma]{Corollary}
\newtheorem{proposition}[lemma]{Proposition}
\newtheorem{theorem}[lemma]{Theorem}

\newenvironment{tenumerate}{
  \begin{enumerate}
  
  }{\end{enumerate}}
\def\Ga{{\mathfrak{a}}}
\def\Gb{{\mathfrak{b}}}
\def\Gg{{\mathfrak{g}}}
\def\Gh{{\mathfrak{h}}}

\def\Gl{{\mathfrak{l}}}

\def\Gn{{\mathfrak{n}}}
\def\Gp{{\mathfrak{p}}}

\def\BC{{\mathbb{C}}}

\def\BM{{\mathbb{M}}}

\def\BO{{\mathbb{O}}}

\def\BQ{{\mathbb{Q}}}
\def\BR{{\mathbb{R}}}
\def\BZ{{\mathbb{Z}}}

\def\CA{{\mathcal A}}
\def\CC{{\mathcal{C}}}
\def\CD{{\mathcal D}}
\def\CO{{\mathcal O}}

\def\lam{\lambda}
\def\lan{\langle}
\def\ran{\rangle}
\def\nn{\nonumber}

\newcommand{\e}{{\rm e}}

\def\cl{{\rm{cl}}}
\def\mod{{\rm{mod}}}

\def\ad{\mathop{\rm ad}\nolimits}

\def\ch{\mathop{\rm ch}\nolimits}

\def\Ext{{\mathop{\rm Ext}\nolimits}}

\def\id{\mathop{\rm id}\nolimits}

\def\Image{\mathop{\rm Im}\nolimits}
\def\Ker{\mathop{\rm Ker\hskip.5pt}\nolimits}
\def\limi{\mathop{\mathop{\rm lim}_{\longrightarrow}}}

\def\Ob{\mathop{\rm Ob}\nolimits}
\def\re{{\rm re}}

\def\Wt{{\rm Wt}}

\newcommand{\U}{{U(\Gg)}}
\newcommand{\iso}{\mathrel{\longrightarrow{\kern-18pt\raise3.5pt\hbox{$\sim$
}}}\enspace}
\def\endeq{\end{eqnarray}}
\def\endeqn{\end{eqnarray*}}
\def\eq{\begin{eqnarray}}
\def\eqn{\begin{eqnarray*}}

\title[Characters of irreducible modules]
{{Characters of irreducible modules with}\\
{non-critical highest weights over affine Lie algebras}}
\author[M.~Kashiwara]{Masaki Kashiwara}
\address[M.~Kashiwara]{Research Institute for Mathematical Sciences, Kyoto
University, Kyoto, 606--8502, Japan}
\author[T.~Tanisaki]{Toshiyuki Tanisaki}
\address[T.~Tanisaki]{Department of Mathematics, Hiroshima University,
Higashi-Hiroshima 739--8526 Japan}
\begin{document}
\begin{abstract}
We shall derive Kazhdan-Lusztig type character formula for the irreducible
modules with arbitrary non-critical highest weights over affine Lie algebras
from the rational case by using the translation functor, the Enright functor
and Jantzen's deformation argument.
\end{abstract}
\maketitle

\section{Introduction}
\setcounter{equation}{0}
The aim of this paper is to give a character formula for the irreducible
modules with arbitrary non-critical highest weights over affine Lie
algebras.

Let us first recall the history of the corresponding problem for
finite-dimensional semisimple Lie algebras.
In \cite{KL1} Kazhdan-Lusztig proposed a conjecture describing the
characters of the irreducible modules with integral highest weights over
finite-dimensional semisimple Lie algebras
in terms of Kazhdan-Lusztig polynomials.
This conjecture was proved by Beilinson-Bernstein~\cite{BB} and
Brylinski-Kashiwara~\cite{BK} independently using $D$-modules on the flag
manifolds.
Later its generalization to rational highest weights was obtained by
combining an unpublished result of Beilinson-Bernstein and a result in
Lusztig~\cite{L0}.
Then by a result of Jantzen~\cite{Jantzen} the character formula of the
irreducible modules with arbitrary highest weights is obtained by reducing it to the rational highest weight case.

As for affine Lie algebras, we know already descriptions of the characters
of the irreducible modules with rational non-critical highest weights by
Kashiwara-Tanisaki~\cite{KTneg2}, \cite{KTpos3} (see
Kashiwara (-Tanisaki)~\cite{Kpos1}, \cite{KTpos2},
Kashiwara-Tanisaki~\cite{KTneg1}, and Casian~\cite{C1}, \cite{C2} for the
integral case).
In this paper we shall derive the character formula
for arbitrary non-critical highest
weights over affine Lie algebras from the rational non-critical case by
using the translation functor, the Enright functor and Jamtzen's deformation argument.

Let us describe our results more precisely.
Let $\Gg$ be a finite-dimensional semisimple or
affine Lie algebra over the complex number field $\BC$ with
Cartan subalgebra $\Gh$.
Let $\{\alpha_i\}_{i\in I}$ be the set of simple roots, and let $W$ be the
Weyl group.
For a real root $\alpha$ we denote by $s_\alpha\in W$ the corresponding
reflection.
Fix a 
$W$-invariant
non-degenerate symmetric bilinear form $(\ ,\ )$ on
$\Gh^*$.
Set $\alpha^\vee=2\alpha/(\alpha,\alpha)$ for a real root $\alpha$.
Fix $\rho\in\Gh^*$ satisfying $(\alpha^\vee_i,\rho)=1$ for any $i\in I$, and
define a shifted action of $W$ on $\Gh^*$ by
$$
w\circ\lam=w(\lam+\rho)-\rho
\quad\mbox{for any $\lam\in\Gh^*$}.
$$
When $\Gg$ is affine,
we denote by $\delta$ the positive imaginary root such that any imaginary
root is an integral multiple of $\delta$.

For $\lam\in\Gh^*$ we denote by $\Delta^+(\lam)$ the set of positive real
roots $\alpha$ satisfying $(\alpha^\vee,\lam+\rho)\in\BZ$, and by
$\Pi(\lam)$ the set of $\alpha\in\Delta^+(\lam)$ satisfying
$s_\alpha(\Delta^+(\lam)\setminus\{\alpha\})=\Delta^+(\lam)\setminus\{\alpha
\}$.
Then the subgroup $W(\lam)$ of $W$ generated by
$\{s_\alpha\,;\,\alpha\in\Delta^+(\lam)\}$ is a Coxeter group with the
canonical generator system $\{s_\alpha\,;\,\alpha\in\Pi(\lam)\}$.
We denote the Bruhat ordering and the length function of $W(\lam)$ by
$\geq_\lam$ and $\ell_\lam:W(\lam)\to\BZ_{\geq0}$ respectively.
For $y, w\in W(\lam)$ we denote by $P^\lam_{y,w}(q)\in\BZ[q]$ the
corresponding Kazhdan-Lusztig polynomial (see Kazhdan-Lusztig~\cite{KL1}),
and by $Q^\lam_{y,w}(q)\in\BZ[q]$ the inverse Kazhdan-Lusztig polynomial
defined by
$$
\sum_{x\leq_\lam y\leq_\lam z}
(-1)^{\ell_\lam(y)-\ell_\lam(x)}
Q^\lam_{x,y}(q)P^\lam_{y,z}(q)
=\delta_{x,z}\qquad\mbox{for any $x,z\in W(\lam)$.}
$$
We denote by $W_0(\lam)$ the subgroup of $W(\lam)$ generated by
$\{s_\alpha\,;\,\alpha\in\Delta^+,\,(\alpha^\vee,\lam+\rho)=0\}$.

For $\lam\in\Gh^*$ let $M(\lam)$  (resp.\ $L(\lam)$) be the Verma module
(resp.\ irreducible module) with highest weight $\lam$.
We denote the characters of $M(\lam)$ and $L(\lam)$ by $\ch(M(\lam))$ and
$\ch(L(\lam))$ respectively.
The aim of this paper is to give a description of $\ch(L(\lam))$ for any
$\lam\in\Gh^*$ (satisfying $(\delta,\lambda+\rho)\ne0$ when $\Gg$ is affine).

Set
\begin{equation*}
\begin{split}
\CC&=
 \left\{
 \begin{array}{ll}
  \Gh^*
  &\mbox{when $\Gg$ is finite-dimensional semisimple,}
 \\
  \{\lam\in\Gh^*\,;\,(\delta,\lambda+\rho)\ne0\}
  &\mbox{when $\Gg$ is affine,}
 \end{array}\right.
\\
\CC^+&=
\{\lam\in\CC\,;\,(\alpha^\vee,\lam+\rho)\geq0
\mbox{ for any $\alpha\in\Delta^+(\lam)$}\},
\\
\CC^-&=\{\lam\in\CC\,;\,(\alpha^\vee,\lam+\rho)\leq0
\mbox{ for any $\alpha\in\Delta^+(\lam)$}\}.
\end{split}
\end{equation*}

Let $\lam\in\CC$.
Then $W_0(\lam)$ is a finite group, and we have
$(W(\lam)\circ\lam)\cap(\CC^+\cup\CC^-)\ne\emptyset$.
(see \S\ref{Section:integral root system} below).
Moreover, for any $w\in W(\lam)$ there exists a unique $x\in wW_0(\lam)$
such that its length $\ell_\lam(x)$ is the largest (resp.\ smallest) among the
elements of $wW_0(\lam)$.
We call it the longest (resp.\ shortest) element of $wW_0(\lam)$.

Our main result is the following.
\begin{theorem}
\label{Main theorem}
Let $\Gg$ be a finite-dimensional semisimple or affine Lie algebra.
\begin{tenumerate}
\item
Let $\lam\in\CC^+$.
For any $w\in W(\lam)$ which is the longest element of $wW_0(\lam)$ we have
$$
\ch(L(w\circ\lam))=
\sum_{W(\lam)\ni y\geq_\lam w}
(-1)^{\ell_\lam(y)-\ell_\lam(w)}
Q^\lam_{w,y}(1)\ch(M(y\circ\lam)).
$$
\item
Let $\lam\in\CC^-$.
For any $w\in W(\lam)$ which is the shortest element of $wW_0(\lam)$ we have
$$
\ch(L(w\circ\lam))=
\sum_{W(\lam)\ni y\leq_\lam w}
(-1)^{\ell_\lam(w)-\ell_\lam(y)}
P^\lam_{y,w}(1)\ch(M(y\circ\lam)).
$$
\end{tenumerate}
\end{theorem}

We would like to thank J. Bernstein and J. C. Jantzen for useful conversation and comments.

\section{Integral root systems}
\label{Section:integral root system}
\setcounter{equation}{0}

Since the finite-dimensional case is similar and simpler,
we assume in the sequel that $\Gg$ is affine.
Let $\Gg$ be an affine Lie algebra over the complex number field $\BC$.
Let $\Gh$ be the Cartan subalgebra, and let $\{\alpha_i\}_{i\in
I}\subset\Gh^*$ and $\{h_i\}_{i\in I}\subset\Gh$ be the set of simple roots
and the set of simple coroots respectively.
We assume that $\{\alpha_i\}_{i\in I}$ and $\{h_i\}_{i\in I}$
are linearly independent and $\dim \Gh=|I|+1$.
We denote by $\Delta$ (resp.\ $\Delta_{\rm re}$, $\Delta_{\rm im}$,
$\Delta^+$, $\Delta^-$) the set of roots (resp.\ real roots, imaginary
roots, positive roots, negative roots).
Set $\Delta^\pm_{\rm re}=\Delta_{\rm re}\cap\Delta^\pm$, 
$\Delta^\pm_{\rm im}=\Delta_{\rm im}\cap\Delta^\pm$.
There exists a unique $\delta\in\Delta^+_{\rm im}$ satisfying 
$\Delta^+_{\rm im}=\BZ_{>0}\delta$.
Let $c\in\sum_{i\in I}\BZ_{>0}h_i$ be the central element of $\Gg$ such that
$\BZ\,c=\{h\in\sum_{i\in I}\BZ h_i\,;
\,\lan h, \alpha_i\ran=0\ \mbox{for any $i\in I$}\}$.
Here, $\langle\ ,\ \rangle:\Gh\times\Gh^*\to\BC$ 
denotes the canonical paring.
We set
\begin{equation}
Q=\sum_{i\in I}\BZ\alpha_i\qquad\mbox{and}\qquad
Q^+=\sum_{i\in I}\BZ_{\geq0}\alpha_i.
\end{equation}

We fix a $\BZ$-lattice $P$ of $\Gh^*$ satisfying
\begin{align}
&\alpha_i\in P,\quad\lan h_i,P\ran\subset\BZ,\\
&\mbox{there exists some $\lam\in P$ such that $\lan
h_j,\lam\ran=\delta_{ij}$ for $j\in I$
}
\end{align}
for any $i\in I$.
Set
\begin{align}
P^+&=\{\lam\in P\,;\,\lan h_i,\lam\ran\geq0
\,\,\mbox{for any $i\in I$}\},\\
\Gh_\BQ^*&=\BQ\otimes_\BZ P\subset\Gh^*,\\
\Gh_\BR^*&=\BR\otimes_\BZ P\subset\Gh^*.
\end{align}
We further fix a non-degenerate symmetric bilinear form
$(\ ,\ ):\Gh_\BQ^*\times\Gh_\BQ^*\to\BQ$ satisfying
\begin{align}
&\langle h_i,\lam\rangle=2(\lam,\alpha_i)/(\alpha_i,\alpha_i)\,\,
\mbox{ for any $i\in I$ and $\lam\in\Gh_\BQ^*$}
\end{align}
normalized by
\begin{equation}
\lan c, \lam\ran=(\delta,\lam)
\ \mbox{ for any $\lam\in\Gh_\BQ^*$.}
\end{equation}
Then we have
\begin{equation}
(\alpha,\alpha)/2=
\mbox{$1/3$, $1/2$, $1$, $2$ or $3\ $ for any
$\alpha\in\Delta_\re$.}
\end{equation}
Its scalar extension to $\Gh^*$ is also denoted by
$(\ ,\ ):\Gh^*\times\Gh^*\to\BC$.

For $\alpha\in\Delta_{\rm re}$ we set
\begin{equation}
\alpha^\vee=2\alpha/(\alpha,\alpha)\in\Gh_\BQ^*,
\end{equation}
and define $s_\alpha\in GL(\Gh^*)$ by
\begin{equation}
s_\alpha(\lam)=\lam-(\alpha^\vee,\lam)\alpha\quad
\mbox{for any $\lam\in\Gh^*$.}
\end{equation}
The subgroup $W$ of $GL(\Gh^*)$ generated by
$\{s_{\alpha}\,;\,\alpha\in\Delta_{\rm re}\}$ is called the Weyl group.
It is a Coxeter group with a canonical generator system
$\{s_{\alpha_i}\,;\,i\in I\}$.
We denote its length function by $\ell:W\to\BZ_{\geq0}$.

Fix $\rho\in P$ satisfying $(\alpha^\vee_i,\rho)=1$ for any $i\in I$, and
define a shifted action of $W$ on $\Gh^*$ by
\begin{equation}
w\circ\lam=w(\lam+\rho)-\rho
\quad\mbox{for any $\lam\in\Gh^*$}.
\end{equation}

For a subset $\Gamma$ of $\Gh^*$ 
we denote by $\BC\Gamma$ (resp.\ $\BR\Gamma$, $\BQ\Gamma$) the vector
subspace of $\Gh^*$ over $\BC$ (resp.\ $\BR$, $\BQ$) 
spanned by $\Gamma$.

Set
\begin{equation}
E=\BR\Delta_{\rm re}=\{\lam\in \Gh_\BR^*\,;\,(\delta,\lam)=0\},\qquad
E_\cl=E/\BR\delta,
\end{equation}
and let $\cl:E\to E_\cl$ denote the projection.
The restriction $(\ ,\ ):E\times E\to\BR$ of
$(\ ,\ ):\Gh_\BR^*\times\Gh_\BR^*\to\BR$ is positive semi-definite with
radical $\BR\delta$.
Thus it induces a positive definite symmetric bilinear form
$(\ ,\ ):E_\cl\times E_\cl\to\BR$.
Set ${\Delta}_\cl=\cl(\Delta_{\rm re})$.
Then ${\Delta}_\cl$ is a (not necessarily reduced) finite root system in
$E_\cl$.

For each $\gamma\in\Delta_\cl$ there exists some 
$\tilde{\gamma}\in\Delta_{\rm re}$ and $r_\gamma\in\BZ_{>0}$ satisfying
\begin{align}
\label{real roots}
\cl^{-1}(\gamma)\cap\Delta_{\rm re}
&=\{\tilde{\gamma}+nr_\gamma\delta\,;\,n\in\BZ\},\\
\label{real positive roots}
\cl^{-1}(\gamma)\cap\Delta_{\rm re}^+
&=\{\tilde{\gamma}+nr_\gamma\delta\,;\,n\in\BZ_{\geq0}\},\\
\label{real negative roots}
\cl^{-1}(\gamma)\cap\Delta_{\rm re}^-
&=\{\tilde{\gamma}+nr_\gamma\delta\,;\,n\in\BZ_{<0}\}.
\end{align}
Thus we have
\begin{align}
\label{real roots2}
\Delta_{\rm re}
&=\{\tilde{\gamma}+nr_\gamma\delta\,;\,\gamma\in\Delta_\cl, n\in\BZ\},\\
\label{real positive roots2}
\Delta_{\rm re}^+
&=\{\tilde{\gamma}+nr_\gamma\delta\,;\,\gamma\in\Delta_\cl,
n\in\BZ_{\geq0}\},\\
\label{real negative roots2}
\Delta_{\rm re}^-
&=\{\tilde{\gamma}+nr_\gamma\delta\,;\,\gamma\in\Delta_\cl, n\in\BZ_{<0}\}.
\end{align}
We have $\BZ r_\gamma=\BZ\cap\BZ(\gamma,\gamma)/2$.

We call a subset $\Delta_1$ of $\Delta_{\rm re}$ a {\em subsystem}
of $\Delta_{\rm re}$
if $s_\alpha\beta\in\Delta_1$ for any $\alpha, \beta\in\Delta_1$ (see
Kashiwara-Tanisaki~\cite{KTpos3} and Moody-Pianzola~\cite{MP}).
For a subsystem $\Delta_1$ of $\Delta_{\rm re}$ we set
\begin{align}
\Delta_1^\pm&=\Delta^\pm\cap\Delta_1,\\
\Pi_1&=\{\alpha\in\Delta_1^+\,;\,
s_\alpha(\Delta_1^+\setminus\{\alpha\})\subset\Delta_1^+\},\\
W_1&=\lan s_\alpha\,;\,\alpha\in\Delta_1\ran,\\
S_1&=\{s_\alpha\,;\,\alpha\in\Pi_1\}.
\end{align}
We call the elements of $\Delta_1^+$ (resp.\ $\Delta_1^-$, $\Pi_1$) positive
roots (resp.\ negative roots, simple roots) for $\Delta_1$, and $W_1$ the
Weyl group for $\Delta_1$.
The group $W_1$ is a Coxeter group with a canonical generator system $S_1$,
and its length function $\ell_1:W_1\to\BZ_{\geq0}$ is given by
$\ell_1(w)=|w\Delta_1^+\cap\Delta_1^-|$.
We have
\begin{equation}
\label{eq:negativity}
\mbox{
$(\alpha,\beta)\leq0$ for any $\alpha, \beta\in\Pi_1$ such that
$\alpha\ne\beta$}
\end{equation}
(see \cite{KTpos3}).
\begin{lemma}
\label{lemma:integral roots4}
The following conditions for a subsystem $\Delta_1$ of $\Delta_{\rm re}$ are
all equivalent to each other.
\begin{tenumerate}
\item
$|\Delta_1|<\infty$,
\item
$|W_1|<\infty$,
\item
$\BC\Delta_1\not\ni\delta$,
\item
$\BQ\Delta_1\not\ni\delta$.
\end{tenumerate}
\end{lemma}
\begin{proof}
It is well-known that (i) and (ii) are equivalent, and
they are also equivalent to the condition that the restriction
$(\ ,\ )\,|\,\BR\Delta_1\times\BR\Delta_1$ of $(\ ,\ ):E\times E\to\BR$ is
positive definite.
Thus the conditions (i) and (ii) are equivalent to $\BR\Delta_1\not\ni\delta$.
This condition is equivalent to (iii) and (iv) because
$\Delta_1\cup\{\delta\}\subset\Gh_\BQ^*\subset\Gh_\BR^*$.
\end{proof}
\begin{lemma}
\label{lemma:positivity of delta}
Let $\Delta_1$ be a subsystem of $\Delta_{\rm re}$ and let $\Pi_1$ be the
set of simple roots for $\Delta_1$.
If $\BQ\Delta_1\ni\delta$, then we have
$\delta\in\sum_{\alpha\in\Pi_1}\BQ_{\geq0}\alpha$.
\end{lemma}
\begin{proof}
Let $\Pi_{2}$ be a minimal subset of $\Pi_{1}$ such that $\BQ\Pi_{2}\ni\delta$.
Write $\delta=\sum_{\alpha\in\Pi_{2}}c_\alpha\alpha$ with $c_\alpha\in\BQ$.
Let $\Pi_{3}=\{\alpha\in\Pi_{2}\,;\,c_\alpha>0\}$, and set
$\gamma=\sum_{\alpha\in\Pi_{3}}c_\alpha\alpha
=\delta+\sum_{\beta\in\Pi_{2}\setminus\Pi_{3}}(-c_\beta)\beta$.
By (\ref{eq:negativity}) we have
$$
0\leq(\gamma,\gamma)=
\sum_{\alpha\in\Pi_{3}}\sum_{\beta\in\Pi_{2}\setminus\Pi_{3}}
c_\alpha(-c_\beta)(\alpha,\beta)\leq0,
$$
and hence $\gamma\in\BQ\delta$.
If $\gamma=0$, then we have
$\delta=\sum_{\beta\in\Pi_{2}\setminus\Pi_{3}}c_\beta\beta
\in
\BQ_{\leq0}\Pi_1
\subset
\sum_{i\in I}\BQ_{\leq0}\alpha_i$.
This is a contradiction.
Thus $\delta\in\BQ\gamma\subset\BQ\Pi_{3}$.
By the minimality of $\Pi_{2}$ we have $\Pi_{2}=\Pi_{3}$, and hence we have
$\delta\in\sum_{\alpha\in\Pi_2}\BQ_{>0}\alpha\subset\sum_{\alpha\in\Pi_1}\BQ
_{\geq0}\alpha$.
\end{proof}
\begin{lemma}
\label{lemma:finiteness of simple roots}
Let $\Pi_1$ be the set of simple roots for a subsystem $\Delta_1$ of
$\Delta_{\rm re}$.
Then we have $|\Pi_1|<\infty$.
\end{lemma}
\begin{proof}
Let $\approx$ be the equivalence relation on $\Pi_1$ generated by
$$
\alpha, \beta\in\Pi_1, (\alpha,\beta)\ne0\Longrightarrow
\alpha\approx\beta,
$$
and let $\{\Pi_{1,a}\,;\,a\in\CA\}$ denote the set of equivalence classes
with respect to $\approx$.

For $a\in\CA$ set $V_a=\BR\Pi_{1,a}$.
Then $\cl(V_a)$ for $a\in\CA$ are all non-zero and mutually orthogonal with
respect to the natural positive definite symmetric bilinear form on $E_\cl$.
Hence $\CA$ is a finite set.
Thus it is sufficient to show that $\Pi_{1,a}$ is a finite set for each
$a\in\CA$.

If $V_a\not\ni\delta$, then $(\ ,\ )|V_a\times V_a$ is positive definite,
and hence $\Delta_\re\cap V_a$ is a finite subsystem of $\Delta_{\rm re}$.
Thus $\Pi_{1,a}$ is a finite set.

Assume that $V_a\ni\delta$.
By Lemma~\ref{lemma:positivity of delta} there exists a finite subset
$\Pi_{2,a}$ of $\Pi_{1,a}$ such that
$\delta=\sum_{\alpha\in\Pi_{2,a}}c_\alpha\alpha$ with
$c_\alpha\in\BQ_{>0}$.
Since
$$
0=(\delta,\beta)
=\sum_{\alpha\in\Pi_{2,a}}c_\alpha(\alpha,\beta)
\qquad
\mbox{for any $\beta\in\Pi_{1,a}\setminus\Pi_{2,a}$,}
$$
(\ref{eq:negativity}) implies
$(\alpha,\beta)=0$ for any $\alpha\in\Pi_{2,a}$ and
$\beta\in\Pi_{1,a}\setminus\Pi_{2,a}$.
Since $\Pi_{1,a}$ is an equivalence class with respect to $\approx$, we
obtain $\Pi_{1,a}=\Pi_{2,a}$.
Therefore, $\Pi_{1,a}$ is a finite set.
\end{proof}

For a subset $J$ of $I$ set
\begin{equation}
\Delta_J=\Delta\cap\sum_{i\in J}\BZ\alpha_i.
\end{equation}
If $J$ is a proper subset of $I$, then $\Delta_J$ is a finite 
subsystem with
$\{\alpha_i\,;\,i\in J\}$ as the set of simple roots.
\begin{lemma}
\label{lemma:integral roots1}
For any finite subsystem $\Delta_1$ of $\Delta$ there exist $w\in W$ and a
proper subset $J$ of $I$ such that  $w\Delta_1\subset\Delta_J$.
\end{lemma}
\begin{proof}
Set $V=\BR\Delta_1$.
By Lemma~\ref{lemma:integral roots4} we have $V\not\ni\delta$.
Since $(\,,\,)|V\times V$ is positive definite, $V\cap\Delta_{\rm \re}$ is a
finite subsystem of $\Delta_{\rm re}$ containing $\Delta_1$.
Hence we can assume $\Delta_1=V\cap\Delta_{\rm \re}$ from the beginning.

Set $V^\perp=\{\mu\in\Gh^*_\BR;(V,\mu)=0\}$. Since $\delta\not\in V$,
$(\delta,\mu)$ is not identically zero on $\mu\in V^\perp$.
Similarly $(\alpha,\mu)$ ($\alpha\in\Delta_\re\setminus \Delta_1$)
is not identically zero on $\mu\in V^\perp$.
Since $\Delta_{\rm re}\setminus\Delta_1$ is a countable set, there exists
some $\lam\in V^\perp$ such that $(\delta,\lam)>0$ and
$(\alpha,\lam)\ne0$ for any $\alpha\in\Delta_{\rm re}\setminus\Delta_1$.
Then we have $\Delta_1=\{\alpha\in\Delta_{\rm re}\,;\,(\alpha,\lam)=0\}$.
Since $(\delta,\lam)>0$, there exist only finitely many
$\alpha\in\Delta_{\rm re}^+$ such that $(\alpha,\lam)<0$ by (\ref{real
positive roots2}).
Hence there exists some $w\in W$ such that $(\alpha,w\lam)\geq0$ for any
$\alpha\in\Delta^+_{\rm re}$ by \cite[Proposition 3.2]{Kac}.
Then we obtain $w\Delta_1=
\{\alpha\in\Delta_{\rm re}\,;\,(\alpha,w\lam)=0\}
=\Delta_J$
with $J=\{i\in I\,;\,(\alpha_i,w\lam)=0\}$.
Since $|\Delta_J|=|\Delta_1|<\infty$, we have $J\ne I$.
\end{proof}

\bigskip
For $\lam\in\Gh^*$ set
\begin{align}
\Delta(\lam)
&=
\{\alpha\in\Delta_{\rm re}\,;\,(\alpha^\vee,\lam+\rho)\in\BZ\},\\
\Delta_0(\lam)
&=
\{\alpha\in\Delta_{\rm re}\,;\,(\alpha^\vee,\lam+\rho)=0\}.
\end{align}
They are subsystems of $\Delta_{\rm re}$.
We denote the set of positive roots, the set of negative roots, the set of
simple roots and the Weyl group for $\Delta(\lam)$ by $\Delta^+(\lam)$,
$\Delta^-(\lam)$, $\Pi(\lam)$ and $W(\lam)$ respectively.
We denote those for $\Delta_0(\lam)$ by $\Delta_0^+(\lam)$,
$\Delta_0^-(\lam)$, $\Pi_0(\lam)$ and $W_0(\lam)$.
The length function for
$W(\lam)$ is denoted by $\ell_\lam:W(\lam)\to\BZ_{\geq0}$.
\begin{lemma}
\label{lemma:integral roots3}
For $\lam\in\Gh^*$ such that $\Delta(\lam)\not=\emptyset$,
the following conditions  are equivalent.
\begin{tenumerate}
\item
$|\Delta(\lam)|<\infty$.
\item
$(\delta,\lam+\rho)\notin\BQ$.
\end{tenumerate}
\end{lemma}
\begin{proof}
(i)$\Rightarrow$(ii).
Assume $(\delta,\lam+\rho)\in\BQ$ and $\Delta(\lam)\ne\emptyset$.
Take $\alpha\in\Delta(\lam)$.
By (\ref{real roots}) there exists some $r\in\BZ_{>0}$ such that $\alpha+\BZ
r\delta\subset\Delta_{\rm re}$.
For $n\in\BZ$ we have
$$
((\alpha+nr\delta)^\vee,\lam+\rho)
=(\alpha^\vee,\lam+\rho)+2nr(\delta,\lam+\rho)/(\alpha,\alpha),
$$
and hence we have $\alpha+nr\delta\in\Delta(\lam)$ for any $n\in\BZ$
satisfying $2nr(\delta,\lam+\rho)/(\alpha,\alpha)\in\BZ$.
Thus $|\Delta(\lam)|=\infty$.

\medskip
\noindent
(ii)$\Rightarrow$(i).
Assume $|\Delta(\lam)|=\infty$.
By Lemma~\ref{lemma:integral roots4} we have $\BQ\Delta(\lam)\ni\delta$.
Then we have
$$
(\delta,\lam+\rho)\in
\sum_{\alpha\in\Delta(\lam)}\BQ(\alpha^\vee,\lam+\rho)
\subset\BQ.
$$
\end{proof}
Set
\begin{equation}
\label{eq:the definition of CC}
\CC=\{\lam\in\Gh^*\,;\,(\delta,\lambda+\rho)\ne0\}.
\end{equation}
\begin{lemma}
For any $\lam\in\CC$ we have $|\Delta_0(\lam)|<\infty$.
\end{lemma}
\begin{proof}
Since $(\delta,\lam+\rho)\ne0$, (\ref{real roots}) implies
$|\cl^{-1}(\gamma)\cap\Delta_0(\lam)|\leq1$ for any $\gamma\in\Delta_\cl$.
Thus we have $|\Delta_0(\lam)|\leq|\Delta_\cl|<\infty$.
\end{proof}
In the sequel, we use the following proposition
on the existence of rational points of a subset
defined by linear inequalities.
Since the proof is elementary, we do not give the proof.

\begin{proposition}
\label{prop:rat.exist}
Let $V_\BQ$ be a finite-dimensional $\BQ$-vector space and set
$V_\BR=\BR\otimes_{\BQ}V_\BQ$ and
$V=\BC\otimes_{\BQ}V_\BQ$.
Let $X$ be a subset of $V_\BQ^*$ and
$\{Y_a\}_{a\in A}$ be a family of non-empty finite subsets
of $V_\BQ^*$. Let $B_x$ $($$x\in X$$)$ and
$C_{y,a}$ $($$a\in A$, $y\in Y_a$$)$ be rational numbers.
Set
\begin{align*}
\Omega&=\{\lam\in V; \lan x,\lam\ran=B_x\quad\mbox{for any $x\in X$}\},\\
\Omega'&=\{\lam\in\Omega\,; \,
\mbox{for any $a\in A$, there exists $y\in Y_a$ such that
$\lan y,\lam\ran\notin C_{y,a}\BZ$}\}.
\end{align*}
\begin{tenumerate}
\item
If $A$ is a finite set and $\Omega'\not=\emptyset$,
then $\Omega'\cap V_\BQ\not=\emptyset$.
\item
If $A$ is a countable set and $\Omega'\not=\emptyset$,
then $\Omega'\cap V_\BR\not=\emptyset$.
Moreover if $z\in V_\BQ^*$ is not contained in the vector subspace $\BQ X$,
then there exists
$\lam\in \Omega'\cap V_\BR$ such that $\lan z,\lam\ran>1$.
\end{tenumerate}
\end{proposition}
\begin{lemma}
\label{lemma:w0}
For any $\lam\in\CC$ we have $W_0(\lam)=\{w\in W\,;\,w\circ\lam=\lam\}$.
\end{lemma}
\begin{proof}
Set $W_1=\{w\in W\,;\,w\circ\lam=\lam\}$.
It is sufficient to show that the group $W_1$ is generated by the
reflections contained in it.
Set
$$
\Omega'=\{\mu\in\CC\,;\,
\mbox{$w\circ\mu=\mu$ for any $w\in W_1$,
$w\circ\mu\ne\mu$ for any $w\in W\setminus W_1$}\}.
$$
Since $\Omega'$ contains $\lam$,
Proposition \ref{prop:rat.exist} (ii)
implies that $\Omega'\cap\Gh_\BR^*$ contains a point $\mu$ such that
$(\delta,\mu+\rho)>0$.
Thus replacing $\lam$ with such a $\mu$,
we may assume that $\lam\in\CC\cap\Gh_\BR^*$ and $(\delta,\lam+\rho)>0$.
Then the assertion follows from \cite[Proposition~3.2]{Kac} and
\cite[Proposition~5.8]{Kac}.
\end{proof}
By a standard argument we have the following.
\begin{lemma}
\label{lemma:dominant chamber2}
Set
\begin{align}
\Gh^{*+}&=\{\lam\in\Gh^*\,;\,(\alpha^\vee,\lam+\rho)\geq0\,\,
\mbox{for any $\alpha\in\Delta^+(\lam)$}
\},\nn\\
\Gh^{*-}&=\{\lam\in\Gh^*\,;\,(\alpha^\vee,\lam+\rho)\leq0\,\,
\mbox{for any $\alpha\in\Delta^+(\lam)$}
\}.\nn
\end{align}
Then for any $\lam\in\Gh^*$, $|(W(\lam)\circ\lam)\cap\Gh^{*\pm}|\le1$.
Moreover, $|(W(\lam)\circ\lam)\cap\Gh^{*+}|=1$ $($resp.\
$|(W(\lam)\circ\lam)\cap\Gh^{*-}|=1)$
if and only if
there exist only finitely many $\alpha\in\Delta^+(\lam)$ satisfying
$(\alpha^\vee,\lam+\rho)<0$ $($resp.\ $(\alpha^\vee,\lam+\rho)>0)$.
\end{lemma}

Set
\begin{align}
\CC^+&=\{\lam\in\CC\,;\,(\alpha^\vee,\lam+\rho)\geq0
\mbox{ for any $\alpha\in\Delta^+(\lam)$}\},\\
\CC^-&=\{\lam\in\CC\,;\,(\alpha^\vee,\lam+\rho)\leq0
\mbox{ for any $\alpha\in\Delta^+(\lam)$}\}.
\end{align}
\begin{lemma}
\label{lemma:integral roots5}
Assume $\lam\in\CC$ satisfies $\Delta(\lam)\ne\emptyset$.
\begin{tenumerate}
\item
If $(\delta,\lam+\rho)\notin\BQ$, then we have
$|(W(\lam)\circ\lam)\cap\CC^+|=|(W(\lam)\circ\lam)\cap\CC^-|=1$.
\item
If $(\delta,\lam+\rho)\in\BQ_{>0}$, then we have
$|(W(\lam)\circ\lam)\cap\CC^+|=1$ and $|(W(\lam)\circ\lam)\cap\CC^-|=0$.
\item
If $(\delta,\lam+\rho)\in\BQ_{<0}$, then we have
$|(W(\lam)\circ\lam)\cap\CC^+|=0$ and $|(W(\lam)\circ\lam)\cap\CC^-|=1$.
\end{tenumerate}
\end{lemma}
\begin{proof}

\noindent
(i)\quad
If $(\delta,\lam+\rho)\notin\BQ$, then we have $|\Delta^+(\lam)|<\infty$ by
Lemma~\ref{lemma:integral roots3}.
Hence we have
$|(W(\lam)\circ\lam)\cap\CC^+|=|(W(\lam)\circ\lam)\cap\CC^-|=1$ by
Lemma~\ref{lemma:dominant chamber2}.

\medskip
\noindent
(ii)\quad Assume $(\delta,\lam+\rho)\in\BQ_{>0}$.
Set
\begin{align*}
\Delta_1&=\{\alpha\in\Delta^+(\lam)\,;\,(\alpha^\vee,\lam+\rho)>0\},\\
\Delta_2&=\{\alpha\in\Delta^+(\lam)\,;\,(\alpha^\vee,\lam+\rho)<0\},\\
\Delta_3&=\{\alpha\in\Delta^+(\lam)\,;\,(\alpha^\vee,\lam+\rho)\leq0\}.
\end{align*}
For each $\gamma\in\Delta_\cl$ there exist only finitely many
$\alpha\in\cl^{-1}(\gamma)\cap\Delta_{\rm re}^+$ satisfying
$(\alpha^\vee,\lam+\rho)\in\BZ_{\leq0}$ by (\ref{real positive roots}).
Since $|\Delta_\cl|<\infty$, we obtain $|\Delta_3|<\infty$.
Thus we have $|\Delta_2|\leq|\Delta_3|<\infty$.
On the other hand we have $|\Delta^+(\lam)|=\infty$ by Lemma
\ref{lemma:integral roots3}, and hence
$|\Delta_1|=|\Delta^+(\lam)\setminus\Delta_3|=\infty$.
Thus we obtain the desired result by Lemma~\ref{lemma:dominant chamber2}.

\medskip
\noindent
The assertion (iii)
follows from (ii) by replacing
$\lam$ with $-\lam-2\rho$.
\end{proof}
\begin{corollary}
\label{corollary:integral roots}
For any $\lam\in\CC$ we have
$(W(\lam)\circ\lam)\cap(\CC^+\cup\CC^-)\ne\emptyset$.
\end{corollary}

\begin{lemma}
\label{lemma:integral roots6}
Let $\lam\in\CC$.
\begin{tenumerate}
\item
If $\BQ\Delta(\lam)\ni\delta$, then there exists some
$\mu\in\CC\cap\Gh^*_\BQ$ such that $(\delta,\mu+\rho)=(\delta,\lam+\rho)$,
$\Delta(\mu)=\Delta(\lam)$ and
$(\alpha^\vee,\mu+\rho)=(\alpha^\vee,\lam+\rho)$ for any
$\alpha\in\Delta(\lam)$.
\item
If $\BQ\Delta(\lam)\not\ni\delta$, then there exists some
$\mu\in\CC\cap\Gh^*_\BR$ such that $(\delta,\mu+\rho)>0$,
$\Delta(\mu)=\Delta(\lam)$ and
$(\alpha^\vee,\mu+\rho)=(\alpha^\vee,\lam+\rho)$ for any
$\alpha\in\Delta(\lam)$.
\end{tenumerate}
\end{lemma}
\begin{proof}
Set
\begin{align}
\Omega&=
\{\mu\in\Gh^*\,;\,
(\alpha^\vee,\mu+\rho)=(\alpha^\vee,\lam+\rho)\,\,
\mbox{ for any $\alpha\in\Delta(\lam)$}\},\nn\\
\Omega'&=
\{\mu\in\Omega\,;\,
(\alpha^\vee,\mu+\rho)\notin\BZ\,\,
\mbox{ for any $\alpha\in\Delta_{\rm re}\setminus\Delta(\lam)$}\}.\nn
\end{align}
Then $\Omega'$ contains $\lam$.

\medskip
\noindent
(i)\quad By the definition of $\Omega$ we have
\begin{equation}
\label{eq:lemma:integral roots6}
\mbox{
$(\gamma,\mu+\rho)=(\gamma,\lam+\rho)$\quad for any $\gamma\in\BQ\Delta(\lam)$
and $\mu\in\Omega$.}
\end{equation}
In particular, we have
\begin{equation}
\label{eq:rationality2}
\mbox{
$(\delta,\mu+\rho)=(\delta,\lam+\rho)\in\BQ$\quad for any $\mu\in\Omega$}
\end{equation}
by $\BQ\Delta(\lam)\ni\delta$.
Thus $\Omega\subset\CC$.
Hence it is sufficient to show $\Omega'\cap\Gh_\BQ^*\ne\emptyset$.

\noindent
Set
$$
\Delta_{\cl,1}
=\{\gamma\in\Delta_\cl\,;\,\cl^{-1}(\gamma)\cap\Delta(\lam)=\emptyset\}.\qquad
\Delta_{\cl,2}=\Delta_{\cl}\setminus\Delta_{\cl,1}.
$$
Let $\mu\in\Omega$.
(\ref{real roots}) and the assumption $\BQ\Delta(\lam)\ni\delta$ imply
$\cl^{-1}(\Delta_{\cl,2})\cap\Delta_{\rm re}\subset\BQ\Delta(\lam)$.
Hence
$(\alpha^\vee,\mu+\rho)=(\alpha^\vee,\lam+\rho)\notin\BZ$ for any
$\alpha\in\cl^{-1}(\Delta_{\cl,2})\cap(\Delta_{\rm re}\setminus\Delta(\lam))$.
Thus we have $\mu\in\Omega'$ if and only if
$(\alpha^\vee,\mu+\rho)\notin\BZ$ for any 
$\alpha\in\Delta_{\rm re}\cap\cl^{-1}(\Delta_{\cl,1})$.
By (\ref{real roots}) and
(\ref{eq:rationality2}), this condition is equivalent to
$$
(\tilde{\gamma}^\vee,\mu+\rho)\notin
\BZ+\frac{2r_\gamma(\delta,\lam+\rho)}{(\gamma,\gamma)}\BZ
\quad
\mbox{for any $\gamma\in\Delta_{\cl,1}$.}
$$
Thus we obtain
\begin{equation*}
\Omega'=
\{\mu\in\Omega\,;\,
(\tilde{\gamma}^\vee,\mu+\rho)\notin
q_{\gamma}\BZ\,\,
\mbox{
for any $\gamma\in\Delta_{\cl,1}$}\},
\end{equation*}
where $\{q_\gamma\,;\,\gamma\in\Delta_{\cl,1}\}$ is
a set of positive rational numbers.
Then $\Omega'$ contains $\lam$,
and Proposition \ref{prop:rat.exist} (i)
implies that $\Omega'\cap\Gh_\BQ^*\ne\emptyset$.

\medskip
\noindent
(ii)\quad
This follows immediately from Proposition \ref{prop:rat.exist} (ii).
\end{proof}
\begin{lemma}
\label{lemma:integral roots2}
For any $\lam\in\CC^+\cup\CC^-$,
there exist $w\in W$ and a proper subset $J$ of $I$
such that $w\Delta^+(\lam)\subset\Delta^+$ and $w\Delta_0(\lam)=\Delta_J$.
\end{lemma}
\begin{proof}
By replacing $\lam$ with $-2\rho-\lam$ if necessary,
we may assume $\lam\in\CC^+$ from the beginning.
Let us first show
that there exists some
$\mu\in\CC\cap\Gh_\BR^*$ such that $(\delta,\mu+\rho)>0$,
$\Delta(\mu)=\Delta(\lam)$ and
$(\alpha^\vee,\mu+\rho)=(\alpha^\vee,\lam+\rho)$ for any
$\alpha\in\Delta(\lam)$.
If $\BQ\Delta(\lam)\ni\delta$, then we have $(\delta,\lam+\rho)>0$ by
Lemma~\ref{lemma:positivity of delta}, and
Lemma~\ref{lemma:integral roots6} (i) implies the existence of such a $\mu$.
If $\BQ\Delta(\lam)\not\ni\delta$,
then Lemma~\ref{lemma:integral roots6} (ii) implies
the existence of such a $\mu$.

By (\ref{real positive roots2}) there exist only finitely many
$\alpha\in\Delta_{\rm re}^+$ such that $(\alpha^\vee,\mu+\rho)<0$.
Thus there exists some $w\in W$ such that
$(\alpha^\vee,w\circ\mu+\rho)\geq0$ for any $\alpha\in\Delta_{\rm re}^+$ by
\cite[Proposition 3.2]{Kac}.
We may assume that $\ell(w)=\min\{\ell(x)\,;\,x\in wW_0(\mu)\}$.
Then we have $w(\Delta_0^+(\mu))\subset\Delta^+$ by \cite[Proposition
2.2.11]{KTpos3}.
For
$\alpha\in\Delta^+(\mu)\setminus\Delta_0(\mu)=\Delta^+(\lam)\setminus\Delta
_0(\lam)$ we have
$$
(w\alpha^\vee,w\circ\mu+\rho)
=(\alpha^\vee,\mu+\rho)
=(\alpha^\vee,\lam+\rho)>0,
$$
and hence $w\alpha\in\Delta^+$.
Thus we obtain $w\Delta^+(\lam)\subset\Delta^+$.
Moreover, we have
$$
w\Delta_0(\lam)=w\Delta_0(\mu)=\Delta_0(w\circ\mu)=\Delta_J
$$
with $J=\{i\in I\,;\,(\alpha_i^\vee,w\circ\mu+\rho)=0\}$.
Then $J$ is a proper subset of $I$ by $|\Delta_0(\lam)|<\infty$.
\end{proof}

\section{Translation functor}
\setcounter{equation}{0}
In this section we shall give some properties of the translation functor
(see also Deodhar-Gabber-Kac~\cite{DGK}, and Kumar~\cite{Kumar2}).

For a Lie algebra $\Ga$ over $\BC$ we denote its enveloping algebra by
$U(\Ga)$ and the category of $\Ga$-modules by $\BM(\Ga)$.

For an $\Gh$-module $M$ and $\mu\in\Gh^*$ we set
\begin{equation}
M_\mu=\{m\in M\,;\,hm=\langle h,\mu\rangle m
\mbox{ for any $h\in\Gh$}\}.
\end{equation}
An element $\mu$ of $\Gh^*$ is called a weight of
$M$ if $M_\mu\ne0$.
For an $\Gh$-module $M$ satisfying
\begin{equation}
\label{eq:weight1}
M=\bigoplus_{\mu\in\Gh^*}M_\mu\qquad
\mbox{with $\dim M_\mu<\infty$ for any $\mu\in\Gh^*$,}
\end{equation}
we define its character $\ch(M)$ by the formal sum
\begin{equation}
\ch(M)=\sum_{\mu\in\Gh^*}\dim M_\mu\, \e^\mu.
\end{equation}
We denote by $\BO$ the full subcategory of $\BM(\Gg)$ consisting of
$M\in\Ob(\BM(\Gg))$ satisfying (\ref{eq:weight1}) and
\begin{equation}
\label{eq:weight3}
\begin{array}{ll}
&\mbox{for any $\xi\in\Gh^*$ there exist only finitely many $\mu\in\xi+Q^+$
such that}\\
&\mbox{$M_\mu\ne0$.}
\end{array}
\end{equation}

For $\alpha\in\Delta$ let $\Gg_\alpha$ denote the root space corresponding
to $\alpha$, and set
\begin{equation}
\Gn^+=\bigoplus_{\alpha\in\Delta^+}\Gg_{\alpha},\qquad
\Gn^-=\bigoplus_{\alpha\in\Delta^-}\Gg_{\alpha}\qquad\mbox{and}
\qquad\Gb=\Gh\oplus\Gn^+\,.
\end{equation}

For $\lam\in\Gh^*$ define a $\Gg$-module $M(\lam)$, called the Verma module
with highest weight $\lam$, by
\begin{equation}
M(\lam)=U(\Gg)\otimes_{U(\Gb)}\BC_\lam
\end{equation}
where
$\BC_\lam=\BC 1_\lam$ is the one-dimensional $\Gb$-module
given by
$h1_\lam=\lam(h)1_\lam$ for $h\in\Gh$ and $\Gn^+1_\lam=0$.
We denote its unique irreducible quotient by $L(\lam)$.

We have
\begin{equation}
\ch(M(\lam))=\frac{\e^\lam}
{\prod_{\alpha\in\Delta^+}(1-\e^{-\alpha})^{\dim \Gg_\alpha}}\,.
\end{equation}
Moreover, $M(\lam)$ and $L(\lam)$ are objects of $\BO$ for any $\lam\in\Gh^*$.
For $M\in\Ob(\BO)$ and $\lam\in\Gh^*$ we denote by $[M:L(\lam)]$ the
multiplicity of $L(\lam)$ in $M$ (see \cite[\S9.6]{Kac}).

The following result due to Kac-Kazhdan~\cite{KK} is fundamental in the
study of
highest weight modules.
\begin{proposition}
\label{prop:KK0}
Let $\lam, \mu\in\Gh^*$.
Then the following conditions are equivalent.
\begin{tenumerate}
\item
The multiplicity $[M(\lam):L(\mu)]$ is non-zero.
\item
There exists an injective homomorphism $M(\mu)\to M(\lam)$.
\item
There exist a sequence of
positive roots $\{\beta_k\}_{k=1}^l$,
a sequence of positive integers $\{n_k\}_{k=1}^l$
and a sequence of weights $\{\lam_k\}_{k=0}^l$
such that
$\lam_0=\lam$, $\lam_l=\mu$ and
$\lam_k=\lam_{k-1}-n_k\beta_k$,
$2(\beta_k,\lam_{k-1}+\rho)=n_k(\beta_k,\beta_k)$
for $k=1,\ldots,l$.
\end{tenumerate}
\end{proposition}

For a subset $\CD$ of $\Gh^*$ we denote by $\BO\{\CD\}$ the full subcategory
of $\BO$ consisting of $M\in\Ob(\BO)$ satisfying $[M:L(\mu)]=0$ for any
$\mu\in\Gh^*\setminus \CD$.
For $\lam\in\CC$ (see (\ref{eq:the definition of CC}) for the notation) we
set $\BO[\lam]=\BO\{W(\lam)\circ\lam\}$.
We have obviously $L(\lam)\in\Ob(\BO[\lam])$ for any $\lam\in\CC$.

By Proposition~\ref{prop:KK0} we have the following.
\begin{proposition}
\label{prop:KK}
For any $\lam\in\CC$ we have $M(\lam)\in\Ob(\BO[\lam])$.
\end{proposition}
Define an equivalence relation $\sim$ on $\CC$ by
\begin{equation}
\lam\sim\mu\;
\Longleftrightarrow
\mu\in W(\lam)\circ\lam.
\end{equation}
By Kumar~\cite{Kumar} we have the following.
\begin{proposition}
Any $M\in\Ob(\BO\{\CC\})$ is uniquely decomposed as
$$
M=\bigoplus_{\lam\in\CC/\sim}M[\lam],\qquad
M[\lam]\in\Ob(\BO[\lam]).
$$
\end{proposition}
For $\lam\in\CC$ let
\begin{equation}
P_\lam:\BO\{\CC\}\to\BO[\lam]
\end{equation}
be the projection functor given by $P_\lam(M)=M[\lam]$.

\begin{lemma}
Let $\lam, \mu\in\CC$, $\nu\in\Gh^*$, $x\in W$ satisfy $\mu-\lam=x\nu$.
Then we have $M\otimes L(\nu)\in\Ob(\BO\{\CC\})$ for any $M\in\Ob(\BO[\lam])$.
\end{lemma}
\begin{proof}
It is easily seen that $M\otimes L(\nu)\in\Ob(\BO)$.
Hence it is sufficient to show that if $L(\xi)$ appears as a subquotient of
$M\otimes L(\nu)$, then we have $(\delta,\xi+\rho)\ne0$.

We may assume that $M=L(w\circ\lam)$ for $w\in W(\lam)$.
The central element $c$ of $\Gg$ acts on $L(\eta)$
via the multiplication of the scalar $\lan
c,\eta\ran=(\delta,\eta)$ for any $\eta\in\Gh^*$.
For $w\in W(\lam)$ we have $(\delta,w\circ\lam)=(\delta,\lam)$ by the
$W$-invariance of $\delta$, and hence $c$ acts on $L(w\circ\lam)$ via the
multiplication of $(\delta,\lam)$.
Therefore
we have $cu=(\delta,\lam+\nu)u$ for any $u\in M\otimes L(\nu)$.
If $L(\xi)$ appears as a subquotient of $M\otimes L(\nu)$, then we have
$(\delta,\xi)=(\delta,\lam+\nu)$, and hence
$$
(\delta,\xi+\rho)=(\delta,\lam+\nu+\rho)
=(\delta,\lam+x\nu+\rho)=(\delta,\mu+\rho)\ne0.
$$
\end{proof}

For $\lam, \mu\in\CC$ satisfying
\begin{equation}
\label{cond:tr1}
\mu-\lam\in WP^+\,,
\end{equation}
we define a functor
\begin{equation}
T^\lam_\mu:\BO[\lam]\to\BO[\mu]
\end{equation}
by $T^\lam_\mu(M)=P_\mu(M\otimes L(\nu))$, where $\nu$ is a unique
element of $P^+$ such that $\mu-\lam\in W\nu$.
It is obviously an exact functor.

The proofs of Lemma~\ref{lemma:chameber geometry},
Proposition~\ref{prop:tanslation of M} and Proposition~\ref{prop:translation
of L} below are similar to those for finite-dimensional semisimple Lie
algebras given in Jantzen~\cite{Jantzen}.
We reproduce it here for the sake of completeness.

\begin{lemma}
\label{lemma:chameber geometry}
Assume that we have either $\lam, \mu\in\CC^+$ or $\lam, \mu\in\CC^-$ and
that $\mu-\lam\in W\nu$ for $\nu\in P^+$.
Denote by $\Gamma$ the set of weights of $L(\nu)$.
Then for any $w\in W(\lam)$ satisfying $w\circ\mu-\lam\in\Gamma$ we have
$w\in W_0(\lam)W_0(\mu)$.
\end{lemma}
\begin{proof}
By the assumption we have $\Delta(\lam)=\Delta(\mu)$ and $W(\lam)=W(\mu)$.
Assume that there exists some $w\in W(\lam)\setminus W_0(\lam)W_0(\mu)$
satisfying $w\circ\mu-\lam\in\Gamma$.
We may assume that its length $\ell_\lam(w)$ is the smallest among such
elements.
Set $\xi=w\circ\mu-\lam\in\Gamma$.

Since $w$ is the shortest element of $wW_0(\mu)$,
\cite[Proposition 2.2.11]{KTpos3} implies
\begin{equation}
\label{eq1:chambergeom}
w\Delta^+_0(\mu)\subset\Delta^+(\lam).
\end{equation}
Since $w$ is the shortest element of $W_0(\lam)w$,
\begin{equation}
\label{eq2:chambergeom}
w^{-1}\Delta^+_0(\lam)\subset\Delta^+(\lam).
\end{equation}

By $w\ne1$ there exists some $\alpha\in\Delta^+(\lam)$ satisfying
$\ell_\lam(s_\alpha w)<\ell_\lam(w)$.
Then we have $w^{-1}\alpha\in\Delta^-(\lam)$.
Hence we have $\alpha\in\Delta^+(\lam)\setminus\Delta^+_0(\lam)$ by
(\ref{eq2:chambergeom}).
If $w^{-1}\alpha\in\Delta_0(\mu)$, then we have
$-w^{-1}\alpha\in\Delta^+_0(\mu)\cap w^{-1}\Delta^-(\lam)$.
This contradicts (\ref{eq1:chambergeom}).
Thus we obtain $w^{-1}\alpha\in\Delta^-(\mu)\setminus\Delta^-_0(\mu)$.
Set
$$
m=(\alpha^\vee,\lam+\rho), \quad
n=-(w^{-1}\alpha^\vee,\mu+\rho)=-(\alpha^\vee,w(\mu+\rho)).
$$
By $\alpha\in\Delta^+(\lam)\setminus\Delta^+_0(\lam)$ and
$w^{-1}\alpha\in\Delta^-(\mu)\setminus\Delta^-_0(\mu)$ we have
$m, n\in\BZ_{>0}$ if $\lam, \mu\in\CC^+$ and $m, n\in\BZ_{<0}$ if $\lam,
\mu\in\CC^-$.
Now we have
\begin{align}
&s_\alpha w\circ\mu-\lam
=s_\alpha w(\mu+\rho)-w(\mu+\rho)+\xi
=\xi+n\alpha,\nn\\
&s_\alpha\xi
=\xi-(\alpha^\vee,\xi)\alpha
=\xi-((\alpha^\vee,w(\mu+\rho))-(\alpha^\vee,\lam+\rho))\alpha
=\xi+(m+n)\alpha.\nn
\end{align}
Since $\xi$ and $s_\alpha\xi=\xi+(m+n)\alpha$ are elements of $\Gamma$, we
have $s_\alpha w\circ\mu-\lam=\xi+n\alpha\in\Gamma$.
By $\ell_\lam(s_\alpha w)<\ell_\lam(w)$ we obtain $s_\alpha w\in
W_0(\lam)W_0(\mu)$ by the minimality of $\ell_\lam(w)$.
Hence we have $s_\alpha w\circ\mu-\lam\in W_0(\lam)(\mu-\lam)\subset W\nu$.
It follows that $\xi+n\alpha$ is an extremal weight of $L(\nu)$.
This contradicts $\xi$, $\xi+(m+n)\alpha\in\Gamma$, and $m, n\in\BZ_{>0}$ or
$m, n\in\BZ_{<0}$.
\end{proof}

\begin{proposition}
\label{prop:tanslation of M}
Let $\lam, \mu\in\CC$ such that $\mu-\lam\in WP^+$ and
$\Delta_0(\lam)\subset\Delta_0(\mu)$.
Assume that we have either $\lam, \mu\in\CC^+$ or $\lam, \mu\in\CC^-$.
Then we have $T^\lam_\mu(M(w\circ\lam))=M(w\circ\mu)$ for any $w\in W(\lam)$.
\end{proposition}
\begin{proof}
Take $x\in W$ and $\nu\in P^+$ such that $\mu-\lam=x\nu$.
Let $\Gamma$ be the set of weights of $L(\nu)$.
Since
\begin{equation}
\label{eq:tensor}
M(w\circ\lam)\otimes L(\nu)
=U(\Gg)\otimes_{U(\Gb)}(\BC_{w\circ\lam}\otimes L(\nu))
=U(\Gn^-)\otimes_\BC(\BC_{w\circ\lam}\otimes L(\nu))\,,
\end{equation}
we have
$$
\ch(M(w\circ\lam)\otimes L(\nu))
=\sum_{\xi\in\Gamma}\dim L(\nu)_\xi\ch(M(w\circ\lam+\xi)).
$$
This implies
\eqn
\ch(T^\lam_\mu(M(w\circ\lam)))
&=&\sum_{\xi\in\Gamma}\dim L(\nu)_\xi\ch(P_\mu(M(w\circ\lam+\xi)))\\
&=&\sum_{\xi\in\Gamma,\,w\circ\lam+\xi\in W(\mu)\circ\mu}
\dim L(\nu)_\xi\ch((M(w\circ\lam+\xi)))\,.
\endeqn
Assume that $w\circ\lam+\xi=y\circ\mu$ for $\xi\in\Gamma$ and $y\in W(\lam)$.
Then we have $w^{-1}y\circ\mu-\lam=w^{-1}\xi\in\Gamma$, and hence
$w^{-1}y\in W_0(\lam)W_0(\mu)=W_0(\mu)$ by Lemma~\ref{lemma:chameber
geometry}.
Thus we have
$$
\xi=w(\mu-\lam)=wx\nu\qquad\mbox{and}
\quad
w\circ\lam+\xi
=w\circ(\lam+x\nu)
=w\circ\mu.
$$
Hence we obtain $\ch(T^\lam_\mu(M(w\circ\lam)))=\ch(M(w\circ\mu))$.
In particular, there exists some $v\in(M(w\circ\lam)\otimes
L(\nu))_{w\circ\mu}\setminus\{0\}$ such that $\Gn^+v=0$.
By (\ref{eq:tensor}),
$M(w\circ\lam)\otimes L(\nu)$ is a free $U(\Gn^-)$-module.
Thus the morphism $U(\Gn^-)\to M(w\circ\lam)\otimes L(\nu)$ given by
$u\mapsto uv$ is injective.
It follows that $T^\lam_\mu(M(w\circ\lam))$ contains $M(w\circ\mu)$ as a
submodule.
Hence we have $T^\lam_\mu(M(w\circ\lam))=M(w\circ\mu)$.
\end{proof}
\begin{corollary}
\label{cor:tanslation of M}
Let $\lam, \mu\in\CC$ such that $\mu-\lam\in WP^+$ and
$\Delta_0(\lam)\subset\Delta_0(\mu)$.
Assume that we have either $\lam, \mu\in\CC^+$ or $\lam, \mu\in\CC^-$.
For $M\in\Ob(\BO[\lam])$ let us write
\begin{equation}
\label{eq:char}
\ch M=\sum_{w\in W(\lam)}a_w\ch(M(w\circ\lam))
\end{equation}
with integers $a_w$.
Then we have
$$
\ch T_\mu^\lam(M)=\sum_{w\in W(\lam)}a_w\ch(M(w\circ\mu))\,.
$$
\end{corollary}
\begin{proof}
If $\lam\in\CC^-$, then $M$ has finite length.
Therefore we can reduce the assertion to the case where $M=M(y\circ\lam)$
with $y\in W(\lam)$. 
Then the assertion follows from the preceding proposition.

Assume now $\lam\in\CC^+$.
It is enough to show
\begin{equation}
\label{ch trans}
\dim(T_\mu^\lam(M))_\xi=\sum_{w\in W(\lam)}a_w\dim(M(w\circ\mu)_\xi)
\end{equation}
for any $\xi\in \Gh^*$.
Let $\Wt(M)$ be the set of weights of $M$.
We set $\Gh^*_N=\{\lam-\sum_{i\in I}n_i\alpha_i;\sum n_i\ge N\}$.
Since $w\circ\lam=\lam$ implies $w\circ\mu=\mu$
by Lemma \ref{lemma:w0},
we may assume $w$ ranges over $W(\lam)/W_0(\lam)$ in (\ref{eq:char}).
If $\Wt(M)\subset \Gh^*_N$ for a sufficiently large $N$, then
$a_w\ne0$ implies that $l_\lam(w)$ is sufficiently large.
Hence the both sides of (\ref{ch trans}) vanish.
Fixing such an $N$ we shall argue by the descending induction on 
$m$ such that $\Wt(M)\setminus\Gh^*_N\subset \Gh^*_m$.
Let $w\circ\lam$ ($w\in W(\lam)$) be a highest weight of $M$.
Then there is an exact sequence
$$
0\to M_1\to M(w\circ\lam)^{\oplus m}\to M\to M_2\to0\,,
$$
where $\Wt(M_k)$ does not contain $w\circ\lam$ ($k=1,2$).
Hence by the induction hypothesis, (\ref{ch trans})
holds for $M_1$. Arguing by the induction
on the cardinality of $\Wt(M)\setminus\Gh^*_N$, (\ref{ch trans}) holds for
$M_2$.
Since
$T_\mu^\lam(M(w\circ\lam))=M(w\circ\mu)$ by the preceding proposition,
(\ref{ch trans}) holds for $M(w\circ\mu)$.
Then (\ref{ch trans}) holds for $M$
because $T_\mu^\lam$ is an exact functor.
\end{proof}
\begin{proposition}
\label{prop:translation of L}
Let $\lam, \mu\in\CC$ such that $\mu-\lam\in WP^+$ and
$\Delta_0(\lam)\subset\Delta_0(\mu)$.
Let $w\in W(\lam)$.
\begin{tenumerate}
\item
If $\lam, \mu\in\CC^+$, then we have
$$
T^\lam_\mu(L(w\circ\lam))=\left\{
\begin{array}{ll}
L(w\circ\mu)&
\mbox{
if $w(\Delta_0^+(\mu)\setminus\Delta_0^+(\lam))\subset\Delta^-(\lam)$,
}\\
0&
\mbox{
otherwise.
}
\end{array}
\right.
$$
\item
If $\lam, \mu\in\CC^-$, then we have
$$
T^\lam_\mu(L(w\circ\lam))=\left\{
\begin{array}{ll}
L(w\circ\mu)&
\mbox{
if $w(\Delta_0^+(\mu)\setminus\Delta_0^+(\lam))\subset\Delta^+(\lam)$,
}\\
0&
\mbox{
otherwise.
}
\end{array}
\right.
$$
\end{tenumerate}
\end{proposition}
\begin{proof}
Since $T^\lam_\mu$ is an exact functor, $T^\lam_\mu(L(w\circ\lam))$ is a
quotient of $T^\lam_\mu(M(w\circ\lam))=M(w\circ\mu)$.
By restricting the non-degenerate contravariant form on
$L(w\circ\lam)\otimes L(\nu)$ we obtain a non-degenerate contravariant form
on $T^\lam_\mu(L(w\circ\lam))$.
Thus we have either $T^\lam_\mu(L(w\circ\lam))=L(w\circ\mu)$ or
$T^\lam_\mu(L(w\circ\lam))=0$.

Assume
$w(\Delta_0^+(\mu)\setminus\Delta_0^+(\lam))\not\subset\Delta^-(\lam)$ in
the case $\lam, \mu\in\CC^+$ and
$w(\Delta_0^+(\mu)\setminus\Delta_0^+(\lam))\not\subset\Delta^+(\lam)$ in
the case $\lam, \mu\in\CC^-$.
Then there exists $\alpha\in\Delta(\lam)$ such that $w\alpha\in
\Delta^+(\lam)$, $(\alpha^\vee,\lam+\rho)>0$, and
$(\alpha^\vee,\mu+\rho)=0$.
Set $\beta=w\alpha\in \Delta^+(\lam)$.
Then we have
$(\beta^\vee,w\circ\lam+\rho)>0$ and $(\beta^\vee,w\circ\mu+\rho)=0$.
By Proposition~\ref{prop:KK0} we have exact sequences
\begin{gather}
0\to M(s_\beta w\circ\lam)\to M(w\circ\lam) \to L\to 0,\nn\\
L\to L(w\circ\lam)\to 0.
\nn
\end{gather}
By applying the exact functor $T^\lam_\mu$,
we obtain exact sequences
\begin{gather}
0\to M(s_\beta w\circ\mu)\to M(w\circ\mu) \to T^\lam_\mu(L)\to 0,
\nn\\
T^\lam_\mu(L)\to T^\lam_\mu(L(w\circ\lam))\to 0.\nn
\end{gather}
Since $M(s_\beta w\circ\mu)\to M(w\circ\mu)$ is an isomorphism, we have
$T^\lam_\mu(L(w\circ\lam))=0$.

Next assume
$w(\Delta_0^+(\mu)\setminus\Delta_0^+(\lam))\subset\Delta^-(\lam)$ in the
case $\lam, \mu\in\CC^+$ and
$w(\Delta_0^+(\mu)\setminus\Delta_0^+(\lam))\subset\Delta^+(\lam)$ in the
case $\lam, \mu\in\CC^-$.
Then we have
\begin{equation}
\label{prperty:w}
\begin{array}{ll}
&\mbox{
$w\alpha\in \Delta^-(\lam)$ for any $\alpha\in\Delta(\lam)$
satisfying $(\alpha^\vee,\lam+\rho)>0$ and}\\
&\mbox{$(\alpha^\vee,\mu+\rho)=0$.}
\end{array}
\end{equation}
Let $M$ be the maximal proper submodule of $M(w\circ\lam)$.
By applying $T^\lam_\mu$ to the exact sequence
$$
0\to M\to M(w\circ\lam)\to L(w\circ\lam)\to0\,,
$$
we obtain an exact sequence
$$
0\to T^\lam_\mu(M)\to M(w\circ\mu)\to T^\lam_\mu(L(w\circ\lam))\to 0.
$$
Thus it is sufficient to show $[T^\lam_\mu(M):L(w\circ\mu)]=0$.
Hence we have only to prove $[T^\lam_\mu(L(z\circ\lam)):L(w\circ\mu)]=0$ for
any $z\in W(\lam)$ satisfying $[M:L(z\circ\lam)]\ne0$.
By Proposition~\ref{prop:KK0} there exists some $\beta\in \Delta^+(\lam)$
such that $(\beta^\vee,w(\lam+\rho))>0$ and $[M(s_\beta
w\circ\lam):L(z\circ\lam)]\ne0$.
For such a $\beta$, $T^\lam_\mu(L(z\circ\lam))$ is a subquotient of
$T^\lam_\mu(M(s_\beta w\circ\lam))=M(s_\beta w\circ\mu)$.
Therefore it is sufficient to show $[M(s_\beta w\circ\mu):L(w\circ\mu)]=0$
for any $\beta\in \Delta^+(\lam)$ such that $(\beta^\vee,w(\lam+\rho))>0$.
Set $\alpha=w^{-1}\beta$.
Then we have $\alpha\in\Delta(\lam)$, $w\alpha\in\Delta^+(\lam)$ and
$(\alpha^\vee,\lam+\rho)>0$.
Since $\alpha\in\Delta^\pm(\lam)$ according to
$\lam$, $\mu\in\CC^\pm$, we have
$(\alpha^\vee,\mu+\rho)\ge0$.
Hence (\ref{prperty:w}) implies
$(\beta^\vee,w(\mu+\rho))=(\alpha^\vee,\mu+\rho)>0$.
Thus we obtain $[M(s_\beta w\circ\mu):L(w\circ\mu)]=0$.
\end{proof}
\begin{proposition}
\label{prop:independence}
Let $\lam_1, \lam_2\in\CC$ such that $\lam_1-\lam_2\in P$ and
$\Delta_0(\lam_1)=\Delta_0(\lam_2)$.
Assume that we have either $\lam_1, \lam_2\in\CC^+$ or $\lam_1,
\lam_2\in\CC^-$.
Let $w\in W(\lam_1)$, and write
$$
\ch(L(w\circ\lam_1))=\sum_{y\in
W(\lam_1)/W_0(\lam_1)}a_{y}\ch(M(y\circ\lam_1))
$$
with $a_{y}\in\BZ$.
Then we have
$$
\ch(L(w\circ\lam_2))=\sum_{y\in
W(\lam_1)/W_0(\lam_1)}a_{y}\ch(M(y\circ\lam_2)).
$$
\end{proposition}
\begin{proof}
Note that $\Delta(\lam_1)=\Delta(\lam_2)$, $W(\lam_1)=W(\lam_2)$
and $W_0(\lam_1)=W_0(\lam_2)$.

\medskip
\noindent
{\sl Case $1$.}\quad $\lam_1, \lam_2\in\CC^+$.

By Lemma~\ref{lemma:integral roots2} there exist $x\in W$ and a proper
subset $J$ of  $I$ such that $x^{-1}\Delta^+(\lam_k)\subset\Delta^+$ and
$x^{-1}\Delta_0(\lam_k)=\Delta_J$ for $k=1$ (and hence also for $k=2$).
Take $\xi_1\in P^+$ such that $(\alpha_i^\vee,\xi_1)=0$ for $i\in J$ and
$(\alpha_i^\vee,\xi_1)\in\BZ_{>0}$ for $i\in I\setminus J$.
Set $\xi_2=\xi_1+x^{-1}(\lam_1-\lam_2)$, $\mu=\lam_1+x\xi_1=\lam_2+x\xi_2$.
Then we have
\begin{align*}
&(\alpha_i^\vee,\xi_2)
=(\alpha_i^\vee,x^{-1}(\lam_1-\lam_2))
=(x\alpha^\vee_i,\lam_1+\rho)-(x\alpha^\vee_i,\lam_2+\rho)=0
\qquad
\mbox{for $i\in J$},\\
&(\alpha_i^\vee,\xi_2)
=(\alpha_i^\vee,\xi_1)+(\alpha_i^\vee,x^{-1}(\lam_1-\lam_2))
\qquad
\mbox{for $i\in I\setminus J$},\\
&(\delta,\mu+\rho)=(\delta,\lam_1+\rho)+\sum_{i\in I}m_i(\alpha_i,\xi_1),
\end{align*}
where $\delta=\sum_{i\in I}m_i\alpha_i$.
By taking $(\alpha_i^\vee,\xi_1)$ for $i\in I\setminus J$ sufficiently large,
we may assume that $\xi_2\in P^+$ and $(\delta,\mu+\rho)\ne0$.
Moreover, we have
$$
(\alpha^\vee,\mu+\rho)=(\alpha^\vee,\lam_1+\rho)+(x^{-1}\alpha^\vee,\xi_1)
\ge0
$$
for any $\alpha\in\Delta^+(\mu)=\Delta^+(\lam_1)$,
and hence we have $\mu\in\CC^+$ and
$\Delta_0(\mu)=\Delta_0(\lam_1)=\Delta_0(\lam_2)$.

Thus
Proposition~\ref{prop:translation of L} implies
$T^{\lam_k}_\mu(L(w\circ\lam_k))=L(w\circ\mu)$ for any $w\in W(\lam_k)$ and
$k=1, 2$.
The assertion then follows from Corollary~\ref{cor:tanslation of M}.

\medskip
\noindent
{\sl Case $2$.}\quad $\lam_1, \lam_2\in\CC^-$.

The proof is similar to the one for the case 1.
Take $x\in W$ and a proper subset $J$ of  $I$ such that
$x^{-1}\Delta^+(\lam_k)\subset\Delta^+$ and
$x^{-1}\Delta_0(\lam_k)=\Delta_J$ for $k=1, 2$.
Take $\xi_1\in P^+$ such that $(\alpha_i^\vee,\xi_1)=0$ for $i\in J$ and
$(\alpha_i^\vee,\xi_1)\in\BZ_{>0}$ for $i\in I\setminus J$.
Set $\xi_2=\xi_1-x^{-1}(\lam_1-\lam_2)$, $\mu=\lam_1-x\xi_1=\lam_2-x\xi_2$.
By taking $(\alpha_i^\vee,\xi_1)$ for $i\in I\setminus J$ sufficiently large,
we have $\mu\in\CC^-$, $\xi_2\in P^+$ and $\Delta_0(\mu)=\Delta_0(\lam_k)$
for $k=1, 2$.
Thus Proposition~\ref{prop:translation of L} implies
$T^\mu_{\lam_k}(L(w\circ\mu))=L(w\circ\lam_k)$ for any $w\in W(\lam_k)$ and
$k=1, 2$.
Hence we obtain the desired result by Corollary~\ref{cor:tanslation of M}.
\end{proof}
\begin{proposition}
\label{reduction theorem}
Assume that $\lam$, $\mu\in\CC^+$ $($resp.\ $\lam,\mu\in\CC^-$$)$
satisfy
\begin{equation}
\label{eq:theorem(i)1}
\mu-\lam\in P,\qquad
\Delta_0(\lam)=\emptyset.
\end{equation}
Assume that $w\in W(\lam)$ is the longest $($resp. shortest$)$ element of
$wW_0(\mu)$.
Write
\begin{equation}
\label{eq:theorem(i)2}
\ch(L(w\circ\lam))=\sum_{y\in W(\lam)}a_{y}\ch(M(y\circ\lam))\qquad
\mbox{with $a_{y}\in\BZ$}.
\end{equation}
Then we have
\begin{equation}
\label{eq:theorem(i)3}
\ch(L(w\circ\mu))=\sum_{y\in W(\lam)}a_{y}\ch(M(y\circ\mu)).
\end{equation}
\end{proposition}

\begin{proof}
Let us prove first the case where $\lam,\mu\in\CC^+$.
We first prove the following statement.
\begin{equation}
\label{statement in theorem}
\begin{array}{l}
\mbox{Let $\nu\in\CC^+$. For any $N\in\BZ_{>0}$ there exists some
$\tilde{\nu}\in\CC^+$ such that}\\
\mbox{$\tilde{\nu}-\nu\in P$, $\Delta_0(\tilde{\nu})=\Delta_0(\nu)$,
$(\alpha^\vee,\tilde{\nu}+\rho)\geq
N$ for any $\alpha\in\Delta^+(\nu)\setminus\Delta_0(\nu)$,}\\
\mbox{and $(\delta,\tilde{\nu}+\rho)-(\delta,{\nu}+\rho)\in\BZ_{\geq N}$.}
\end{array}
\end{equation}
By Lemma~\ref{lemma:integral roots2} there exist $x\in W$ and a proper
subset $J$ of $I$ such that $x\Delta^+(\nu)\subset\Delta^+$ and
$x\Delta_0(\nu)=\Delta_J$.
Take $\xi\in P^+$ such that $(\alpha_i^\vee,\xi)=0$ for $i\in J$ and
$(\alpha_i^\vee,\xi)>0$ for $i\in I\setminus J$.
Set $\tilde{\nu}=\nu+x^{-1}\xi$.
Then we have
$(\alpha^\vee,\tilde{\nu}+\rho)=(\alpha^\vee,\nu+\rho)+(x\alpha^\vee,\xi)$
for any $\alpha\in\Delta(\lam)$ and
$(\delta,\tilde{\nu}+\rho)=(\delta,\nu+\rho)+(\delta,\xi)$.
Hence by taking $(\alpha_i^\vee,\xi)>0$ for $i\in I\setminus J$ sufficiently
large, we obtain (\ref{statement in theorem}).

Assume that $\mu\in\CC^+$.
Let $N\in\BZ_{>0}$.
By (\ref{statement in theorem}) there exists $\tilde{\mu}\in\CC^+$ such that
$\tilde{\mu}-\mu\in P$, $\Delta_0(\tilde{\mu})=\Delta_0(\mu)$,
$(\alpha^\vee,\tilde{\mu}+\rho)\geq N$ for any
$\alpha\in\Delta^+(\mu)\setminus\Delta_0(\mu)$, and
$(\delta,\tilde{\mu}+\rho)-(\delta,{\mu}+\rho)\in\BZ_{\geq N}$.
By Lemma~\ref{lemma:integral roots2} there exist $x\in W$ and a proper
subset $J$ of $I$ such that $x\Delta^+(\tilde{\mu})=x\Delta^+({\mu})
\subset\Delta^+$ and
$x\Delta_0(\tilde{\mu})=x\Delta_0({\mu})=\Delta_J$.
Let $w_J$ be the longest element of $W_J$.
Take $\nu\in P^+$ such that $(\alpha_j^\vee,\nu)>0$ for any $j\in J$, and
set $\tilde{\lam}=\tilde{\mu}-x^{-1}w_J\nu$.
Then we have
\begin{equation}
\label{eq:theorem:1}
\tilde{\mu}-\tilde{\lam}\in WP^+.
\end{equation}
Since $(\delta,\tilde{\lam}+\rho)=(\delta,\tilde{\mu}+\rho)-(\delta,\nu)$,
 we have
\begin{equation}
\label{eq:theorem:3}
(\delta,\tilde{\lam}+\rho)\ne0
\end{equation}
when $N$ is sufficiently large.
For any $\alpha\in\Delta^+(\tilde{\mu})$ we have
$(\alpha^\vee,\tilde{\lam}+\rho)=(\alpha^\vee,\tilde{\mu}+\rho)-(w_Jx\alpha^
\vee,\nu)$.
If $\alpha\in\Delta_0^+(\tilde{\mu})=\Delta_0^+({\mu})$, then we have
$(\alpha^\vee,\tilde{\mu}+\rho)=0$ and $w_Jx\alpha\in-\Delta^+_J$, and hence
$(\alpha^\vee,\tilde{\lam}+\rho)
\in\BZ_{>0}$.
If $\alpha\in\Delta^+(\tilde{\mu})\setminus\Delta_0^+(\tilde{\mu})$, then we
have $(\alpha^\vee,\tilde{\lam}+\rho)\in\BZ_{>0}$ when $N$ is sufficiently
large.
Since $\Pi(\tilde{\mu})=\Pi({\mu})$ is a finite set, we have
$(\alpha^\vee,\tilde{\lam}+\rho)>0$ for any $\alpha\in\Pi(\tilde{\mu})$ for
a sufficiently large $N$.
By
$\Delta^+(\tilde{\mu})
\subset\sum_{\alpha\in\Pi(\tilde{\mu})}\BZ_{\geq0}\alpha$
we have
\begin{equation}
\label{eq:theorem:4}
(\alpha^\vee,\tilde{\lam}+\rho)\in\BZ_{>0}\quad
\mbox{for any $\alpha\in\Delta^+(\tilde{\mu})=\Delta^+(\tilde\lam)$}
\end{equation}
when $N$ is sufficiently large.

Take $N$ satisfying (\ref{eq:theorem:3}), (\ref{eq:theorem:4}).
Then we have $\tilde{\lam}\in\CC^+$ and $\tilde{\lam}$ satisfies the
condition (\ref{eq:theorem(i)1}) for $\lam$.
By Proposition~\ref{prop:independence} the integers $a_{y}$ in
(\ref{eq:theorem(i)2}) do not depend on the choice of $\lambda$.
Hence (\ref{eq:theorem(i)2}) holds for $\tilde{\lam}$.
Since $w$ is the longest element of
$wW_0(\mu)=wW_0(\tilde\mu)$ we have
$w\Delta_0^+(\tilde\mu)\subset \Delta^-$, and
Proposition~\ref{prop:translation of L}
implies $T^{\tilde\lam}_{\tilde{\mu}}(L(w\circ\tilde\lam))
=L(w\circ{\tilde{\mu}})$.
Then
Corollary~\ref{cor:tanslation of M} implies
$$
\ch(L(w\circ\tilde{\mu}))=\sum_{y\in W(\lam)}a_{y}\ch(M(y\circ\tilde{\mu})).
$$
The desired result follows then from Proposition~\ref{prop:independence}.

\medskip
As the assertion in the case $\mu\in\CC^-$ is proved similarly,
we shall only give a sketch.
By Proposition~\ref{prop:independence} and an analogue of (\ref{statement in
theorem}) we may assume that $(\alpha^\vee,\mu+\rho)$ for
$\alpha\in\Delta^+(\mu)\setminus\Delta_0(\mu)$ and $(\delta,\mu+\rho)$ are
sufficiently small.
Take $x\in W$ and a proper subset $J$ of $I$ satisfying
$x\Delta^+({\mu})\subset\Delta^+$ and $x\Delta_0({\mu})=\Delta_J$.
Take $\nu\in P^+$ such that $(\alpha_j^\vee,\nu)>0$ for any $j\in J$, and
set $\tilde{\lam}={\mu}-x^{-1}\nu$.
Then we have $\tilde{\lam}\in\CC^-$ and $\tilde{\lam}$ satisfies the
condition (\ref{eq:theorem(i)1}) for $\lam$.
Hence we can take $\tilde{\lam}$ as $\lam$ by
Proposition~\ref{prop:independence}.
Then we have $T^\lam_{{\mu}}(L(w\circ\lam))=L(w\circ{{\mu}})$
by Proposition~\ref{prop:translation of L}.
Hence we obtain the desired result
by Corollary~\ref{cor:tanslation of M}.
\end{proof}

\section{Enright functor}
\setcounter{equation}{0}

We recall certain properties of the Enright functor which will be used later
(see Enright~\cite{Enright}, Deodhar~\cite{Deodhar},
Kashiwara-Tanisaki~\cite[\S 2.4]{KTpos3}).

For $i\in I$ define a subalgebra $\Gg_i$ of $\Gg$ by
$\Gg_i=\Gh\oplus\Gg_{\alpha_i}\oplus\Gg_{-\alpha_i}$.
Take $e_i\in \Gg_{\alpha_i}$, $f_i\in\Gg_{-\alpha_i}$ such that $[e_i,
f_i]=h_i$.
For $a\in\BC$ we denote by $\BM(\Gg_i,a)$ the full subcategory of
$\BM(\Gg_i)$ consisting of $M\in\Ob(\BM(\Gg_i))$ satisfying
\begin{align}
&M=\bigoplus_{\mu\in\Gh^*}M_\mu,\\
&\dim M_\mu=0\,\,\mbox{unless $\lan h_i,\mu\ran\equiv a\ \mod\,\BZ$,}\\
&\dim\BC[e_i]m<\infty\,\,\mbox{ for any $m\in M$}.
\end{align}
For $\mu\in\Gh^*$ let $M_i(\mu)$ be the Verma module for $\Gg_i$ with
highest weight $\mu$.
We fix a highest weight vector $m_\mu$ of $M_i(\mu)$.
\begin{lemma}
\label{lemma:g_i}
Assume $a\notin\BZ$.
For $M\in\Ob(\BM(\Gg_i,a))$ set $N=\bigoplus_{\mu\in\Gh^*}M_\mu^{e_i}\otimes
M_i(\mu)$, where
$$
M_\mu^{e_i}=
\{m\in M_\mu\,;\,e_im=0\}.
$$
Define a linear map $\varphi:N\to M$ by
$$
\varphi(m\otimes f_i^km_\mu)=f_i^km\qquad
\mbox{for $m\in M_\mu^{e_i}$ and $k\in\BZ_{\geq0}$}.
$$
Then $\varphi$ is an isomorphism of $\Gg_i$-modules.
\end{lemma}
\begin{proof}
By the definition of the Verma module $\varphi$ is obviously a homomorphism
of $\Gg_i$-modules.

Let us show that $\varphi$ is surjective.
It is sufficient to show that $M_\xi\subset\Image(\varphi)$ for any
$\xi\in\Gh^*$.
Let $m\in M_\xi$ satisfying $e_i^nm=0$.
We show by induction on $n$ that
$m\in\sum_{k=0}^\infty f_i^kM_{\xi+k\alpha_i}^{e_i}$.
The case $n=0$ is trivial.
Assume $n>0$.
Since $e_i^{n-1}(e_im)=0$, we have $e_im\in\sum_{k=0}^\infty
f_i^kM_{\xi+(k+1)\alpha_i}^{e_i}$ by the hypothesis of induction.
By $a\notin\BZ$ the linear map
$f_i^{k+1}M_{\xi+(k+1)\alpha_i}^{e_i}\to
f_i^kM_{\xi+(k+1)\alpha_i}^{e_i}\,\,(n\mapsto e_in)$ is bijective.
Hence there exists some $u\in\sum_{k=0}^\infty
f_i^{k+1}M_{\xi+(k+1)\alpha_i}^{e_i}$ such that $e_iu=e_im$.
Then we have
$$
m=(m-u)+u\in M_{\xi}^{e_i}+\sum_{k=0}^\infty
f_i^{k+1}M_{\xi+(k+1)\alpha_i}^{e_i}
=\sum_{k=0}^\infty f_i^{k}M_{\xi+k\alpha_i}^{e_i}.
$$

Next let us show that $\varphi$ is injective.
Assume $\Ker(\varphi)\ne0$.
By $a\notin\BZ$ the Verma module $M_i(\mu)$ is irreducible unless
$M_\mu^{e_i}=0$.
Thus there exist subspaces $N(\mu)$ of $M_\mu^{e_i}$ for $\mu\in\Gh^*$ such
that $\Ker(\varphi)=\bigoplus_{\mu\in\Gh^*}N(\mu)\otimes M_i(\mu)$.
Hence there exists some $m\in M_\mu^{e_i}\setminus\{0\}$ such that $m\otimes
M_i(\mu)\subset\Ker(\varphi)$.
Then we have $m=\varphi(m\otimes m_\mu)=0$.
This is a contradiction.
Thus we have $\Ker(\varphi)=0$.
\end{proof}

We denote by $F:\BM(\Gg)\to\BM(\Gg_i)$ the forgetful functor.
For $a\in\BC$ let $\BM_i(\Gg,a)$ be the full subcategory of $\BM(\Gg)$
consisting of $M\in\Ob(\BM(\Gg))$ satisfying $F(M)\in\Ob(\BM(\Gg_i,a))$.

For $a\in\BC$ define a left $U(\Gg)$-module $U(\Gg)f_i^{a+\BZ}$ by
\begin{equation}
U(\Gg)f_i^{a+\BZ}=\limi_nU(\Gg)f_i^{a-n},
\end{equation}
where $U(\Gg)f_i^{a-n}$ is a rank one free $U(\Gg)$-module generated by the
element $f_i^{a-n}$ and the homomorphism $U(\Gg)f_i^{a-n}\to
U(\Gg)f_i^{a-n-1}$ is given by $f_i^{a-n}\mapsto f_if_i^{a-n-1}$.
Then we have a natural $U(\Gg)$-bimodule structure on $U(\Gg)f_i^{a+\BZ}$
whose  right $U(\Gg)$-module structure is given by
\begin{equation}
\label{eq:Enright1}
\begin{split}
f_i^{a+m}P=\sum_{k=0}^{\infty}\genfrac{(}{)}{0pt}{}{a+m}{k}
(\ad(f_i)^kP)&f_i^{a+m-k}\\
&\mbox{for any $m\in\BZ$ and any $P\in\U$}.
\end{split}
\end{equation}
Note that the $U(\Gg)$-bimodule $U(\Gg)f_i^{a+\BZ}$ depends only on
$(a\mbox{ mod } \BZ)\in\BC/\BZ$.

For $M\in\Ob(\BM_i(\Gg,a))$ we set
\begin{equation}
S_i(a)(M)=
\{
m\in U(\Gg)f_i^{a+\BZ}\otimes_{U(\Gg)}M\,;\,
\dim \BC[e_i]m<\infty
\}.
\end{equation}
It defines a left exact functor
\begin{equation}
S_i(a):\BM_i(\Gg,a)\to\BM_i(\Gg,-a),
\end{equation}
called the Enright functor corresponding to $i$.

By the morphism of $U(\Gg)$-bimodules
\begin{equation}
U(\Gg)\to U(\Gg)f_i^{-a+\BZ}\otimes_{U(\Gg)}U(\Gg)f_i^{a+\BZ}\qquad(1\mapsto
f_i^{-a}\otimes f_i^a)
\end{equation}
we obtain a canonical morphism of functors (see \cite[\S 2.4]{KTpos3})
\begin{equation}
\label{eq:T^2}
\id_{\BM_i(\Gg,a)}\to S_i(-a)\circ S_i(a).
\end{equation}

By \cite[\S 2.4]{KTpos3} we have the following result.
\begin{proposition}
\label{prop:Enright:Mw}
Let $\lam\in\Gh^*$, and set $a=\lan h_i,\lam\ran$.
\begin{tenumerate}
\item
If $a\notin\BZ_{>0}$, then we have $S_i(a)(M(\lam))\simeq M(s_i\circ\lam)$.
\item
If $a\notin\BZ$, then the canonical morphism $M(\lam)\to S_i(-a)\circ
S_i(a)(M(\lam))$ induced by $(\ref{eq:T^2})$ is an isomorphism.
\end{tenumerate}
\end{proposition}

We can similarly define a $U(\Gg_i)$-bimodule $U(\Gg_i)f_i^{a+\BZ}$, and the
Enright functor $\overline{S}(a):\BM(\Gg_i,a)\to\BM(\Gg_i,-a)$ for $\Gg_i$
is given by
\begin{equation}
\overline{S}(a)(M)=
\{
m\in U(\Gg_i)f_i^{a+\BZ}\otimes_{U(\Gg_i)}M\,;\,
\dim \BC[e_i]m<\infty
\}
\nn
\end{equation}
for any $M\in\Ob(\BM(\Gg_i,a))$.
Then we have $F\circ {S}_i(a)=\overline{S}(a)\circ F$ by
$U(\Gg_i)f_i^{a+\BZ}\otimes_{U(\Gg_i)}U(\Gg)\simeq U(\Gg)f_i^{a+\BZ}$.

\begin{proposition}
\label{exactness of Enright}
Assume that $a\notin\BZ$.
\begin{tenumerate}
\item
The functor $S_i(a):\BM_i(\Gg,a)\to\BM_i(\Gg,-a)$ gives an equivalence of
categories, and its inverse is given by $S_i(-a)$.
\item
For $\lam\in\Gh^*$ such that $\lan h_i,\lam\ran\equiv a\ \mod\,\BZ$,
we have
$$
S_i(a)(M(\lam))\simeq M(s_i\circ\lam),\qquad
S_i(a)(L(\lam))\simeq L(s_i\circ\lam).
$$
\end{tenumerate}
\end{proposition}
\begin{proof}
(i) We have to show that the canonical morphisms $\id_{\BM_i(\Gg,a)}\to
S_i(-a)\circ S_i(a)$ and $\id_{\BM_i(\Gg,-a)}\to S_i(a)\circ S_i(-a)$ are
isomorphisms.
By the symmetry we have only to show that $\id_{\BM_i(\Gg,a)}\to
S_i(-a)\circ S_i(a)$ is an isomorphism.
Let us show that the canonical morphism $M\to S_i(-a)\circ S_i(a)(M)$ is
bijective for any $M\in\Ob(\BM_i(\Gg,a))$.
By $F\circ S_i(-a)\circ S_i(a)(M)=\overline{S}(-a)\circ\overline{S}(a)\circ
F(M)$ it is sufficient to show that the canonical morphism $N\to
\overline{S}(-a)\circ\overline{S}(a)(N)$ is bijective for any
$N\in\Ob(\BM(\Gg_i,a))$.
This follows from Proposition~\ref{prop:Enright:Mw} for $\Gg_i$ and
Lemma~\ref{lemma:g_i}.

(ii) We have $S_i(a)(M(\lam))\simeq M(s_i\circ\lam)$ by
Proposition~\ref{prop:Enright:Mw}.
By (i) $S_i(a)(L(\lam))$ is the unique irreducible quotient of
$S_i(a)(M(\lam))\simeq M(s_i\circ\lam)$.
Thus we have $S_i(a)(L(\lam))\simeq L(s_i\circ\lam)$.
\end{proof}

\section{Proof of main theorem}
\setcounter{equation}{0}
In this section we shall give a proof of Theorem~\ref{Main theorem}.
We shall use different arguments according to whether
$\BQ\Delta(\lam)\ni\delta$ or not.
Assume $\lam\in\CC^+\cup\CC^-$.
\bigskip

\noindent
{\sl Case $1$.} \quad
$\BQ\Delta(\lam)\ni\delta$.

In this case the following argument is completely similar to Bernstein's
proof of the corresponding result for finite-dimensional semisimple Lie
algebras.

Set
\begin{align}
&\Omega(\lam)=
\{\mu\in\Gh^*\,;\,
(\alpha^\vee,\mu)=(\alpha^\vee,\lam)\quad
\mbox{ for any $\alpha\in\Delta(\lam)$\},}\\
&\Omega'(\lam)=
\{\mu\in\Omega(\lam)\,;\,
(\alpha^\vee,\mu)\notin\BZ\quad
\mbox{ for any $\alpha\in\Delta_{\rm re}\setminus\Delta(\lam)$}\}.
\end{align}
Then we have
\begin{align}
&
\mbox{
$W(\mu)\supset W(\lam)$ and $W_0(\mu)\supset W_0(\lam)$ for any
$\mu\in\Omega(\lam)$,
}\label{pro:omega}\\
&\label{pro:omega2}
\mbox{
$W(\mu)=W(\lam)$ and $W_0(\mu)=W_0(\lam)$ for any $\mu\in\Omega'(\lam)$,
}\\
&
\mbox{
$w\circ\mu-y\circ\mu=w\circ\lam-y\circ\lam$ for any $\mu\in\Omega(\lam)$,
$w,y\in W(\lam)$,
}\\
&
\mbox{$(\delta,\mu)=(\delta,\lam)$ for any $\mu\in\Omega(\lam)$.}
\end{align}
For any $\mu\in\Omega'(\lam)$ and $w\in W(\lam)/W_0(\lam)$ we can write
uniquely
\begin{equation}
\label{eq:awy}
\ch(L(w\circ\mu))=
\sum_{w\in W(\lam)/W_0(\lam)}a_{w,y}(\mu)\ch(M(y\circ\mu))\qquad
\mbox{with $a_{w,y}(\mu)\in\BZ$}
\end{equation}
by Proposition~\ref{prop:KK0} and (\ref{pro:omega2}).

The proof of Proposition~\ref{prop:Jantzen} below is similar to those for finite-dimensional semisimple Lie algebras given in Jantzen~\cite[Theorem 4.9 and Corollar 4.11]{Jantzen}.
We reproduce it here for the sake of completeness.

\begin{proposition}
\label{prop:Jantzen}
For any $w, y\in W(\lam)/W_0(\lam)$ the function $a_{w,y}(\mu)$ defined in
$(\ref{eq:awy})$ is a constant function on $\Omega'(\lam)$.
\end{proposition}
\begin{proof}
For $\mu\in\Omega'(\lam)$ and $w\in W(\lam)/W_0(\lam)$ we have
\begin{equation*}
\begin{split}
\ch(L(w\circ\mu))\,\e^{-w\circ\mu}
&=\sum_{w\in
W(\lam)/W_0(\lam)}a_{w,y}(\mu)\ch(M(y\circ\mu))\,\e^{-w\circ\mu}\\
&=\sum_{w\in
W(\lam)/W_0(\lam)}a_{w,y}(\mu)\e^{y\circ\mu-w\circ\mu}\ch(M(0))\\
&=(\sum_{w\in
W(\lam)/W_0(\lam)}a_{w,y}(\mu)\e^{y\circ\lam-w\circ\lam})\ch(M(0)).
\end{split}
\end{equation*}
Thus for $w\in W(\lam)/W_0(\lam)$ and $\mu, \mu'\in\Omega'(\lam)$ we have
$a_{w,y}(\mu)=a_{w,y}(\mu')$ for any $y\in W(\lam)/W_0(\lam)$
if and only if
$\ch(L(w\circ\mu))\e^{-w\circ\mu}=\ch(L(w\circ\mu'))\e^{-w\circ\mu'}$.
The last condition is equivalent to $\dim L(w\circ\mu)_{w\circ\mu-\xi}=\dim
L(w\circ\mu')_{w\circ\mu'-\xi}$ for any $\xi\in Q^+$.
Fix $w\in W(\lam)/W_0(\lam)$ and $\xi\in Q^+$, and consider the function
\begin{equation}
F(\mu)=\dim L(w\circ\mu)_{w\circ\mu-\xi}
\end{equation}
on $\Omega(\lam)$.
We have only to show that $F$ is constant on $\Omega'(\lam)$.

By a consideration on the contravariant forms on Verma modules we see that
$F$ is a constructible function on $\Omega(\lam)$.
In particular, it is constant on a non-empty Zariski open subset $U$ of
$\Omega(\lam)$.
Let $m$ be the value of $F$ on $U$.
We have to show $F(\mu)=m$ for any $\mu\in\Omega'(\lam)$.
Let $\mu\in\Omega'(\lam)$.
By Proposition~\ref{prop:independence} $a_{w,y}$ is a constant function on
$$
Z=\{\mu'\in\Omega'(\lam)\,;\,\mu'-\mu\in P\}
$$
for any $y\in W(\lam)/W_0(\lam)$.
Thus we see by the above argument that $F$ is constant on $Z$.
Assume for the moment that
\begin{equation}
\label{eq:in Bernstein}
\mbox{
$Z$ is a Zariski dense subset of $\Omega(\lam)$.
}
\end{equation}
Since $Z\cap U\ne\emptyset$, we have $F(\mu')=m$ for some $\mu'\in Z$.
Since $F$ is a constant function on $Z$, we have $F(\nu)=m$ for any $\nu\in Z$.
In particular, we obtain $F(\mu)=m$.

It remains to show (\ref{eq:in Bernstein}).
Set
\begin{align*}
&V=\{\xi\in\Gh^*\,;\,(\alpha^\vee,\xi)=0
\,
\mbox{
for any $\alpha\in\Delta(\lam)$}
\},\\
&V_\BQ=\Gh_\BQ^*\cap V,\\
&V_\BZ=P\cap V.
\end{align*}
We have $\Omega(\lam)=\mu+V$ and $Z=\mu+V_\BZ$.
By the definition of $V$ the natural morphism $\BC\otimes_\BQ V_\BQ\to V$ is
an isomorphism. Since $V_\BQ$ is a $\BQ$-subspace of
$\Gh^*_\BQ=\BQ\otimes_\BZ P$ we have $V_\BQ\simeq\BQ\otimes_\BZ V_\BZ$.
Hence $V_\BZ$ is a $\BZ$-lattice of $V$.
It follows that $Z=\mu+V_\BZ$ is  a Zariski dense subset of
$\Omega(\lam)=\mu+V$.
\end{proof}
Theorem~\ref{Main theorem} is already known to hold
for $\lam\in\Gh^*_\BQ$ such that $\Delta_0(\lam)=\emptyset$
and $\{w\circ\lam=\lam\}=\{1\}$
by Kashiwara-Tanisaki~\cite{KTneg2},
\cite{KTpos3}, and hence for any $\lam\in\Gh^*_\BQ$ by
Lemma~\ref{lemma:w0} and
Proposition~\ref{reduction theorem}.
On the other hand, $\Omega'(\lam)\cap\Gh_\BQ^*\ne\emptyset$ by
Lemma~\ref{lemma:integral roots6}.
Thus the proof of Theorem~\ref{Main theorem} is completed in the case
$\BQ\Delta(\lam)\ni\delta$ by virtue of Proposition~\ref{prop:Jantzen}.

\bigskip

\noindent
{\sl Case $2$.}\quad
$\BQ\Delta(\lam)\not\ni\delta$.

By Lemma~\ref{lemma:finiteness of simple roots}
$\Delta(\lam)$ is a finite set.
Thus by Lemma~\ref{lemma:integral roots1} there exist $x\in W$ and a proper
subset $J$ of $I$ such that $x\Delta(\lam)\subset\Delta_J$.
We may assume that its length $\ell(x)$ is the smallest among the elements
$z\in
W$ satisfying $z\Delta(\lam)\subset\Delta_J$.
Choose a reduced expression $x=s_{\alpha_{i_1}}\cdots s_{\alpha_{i_r}}$ of $x$.
Then we have
\begin{equation}
(\alpha_{i_k}^\vee,s_{\alpha_{i_{k+1}}}\cdots s_{\alpha_{i_r}}\circ\lam+\rho)
\notin\BZ
\,\,\mbox
{ for any $k=1, \dots, r$}.
\end{equation}
Indeed, if $(\alpha_{i_k}^\vee,s_{\alpha_{i_{k+1}}}\cdots
s_{\alpha_{i_r}}\circ\lam+\rho)\in\BZ$, then we have
$\beta=s_{\alpha_{i_{r}}}\cdots
s_{\alpha_{i_{k+1}}}\alpha_{i_k}\in\Delta(\lam)$, and hence
$$
x\Delta(\lam)=
xs_\beta\Delta(\lam)
=
s_{\alpha_{i_{1}}}\cdots s_{\alpha_{i_{k-1}}}s_{\alpha_{i_{k+1}}}\cdots
s_{\alpha_{i_{r}}}\Delta(\lam).
$$
This contradicts the minimality of $\ell(x)$.

Set $\lam'=x\circ\lam$.
Then we have $x\Delta(\lam)=\Delta(\lam')$,
$x\Delta_0(\lam)=\Delta_0(\lam')$ by the definition, and
$x\Pi(\lam)=\Pi(\lam')$ by \cite[Lemma 2.2.2]{KTpos3}.
In particular, $w\mapsto xwx^{-1}$ induces
an isomorphism $W(\lam)\to W(\lam')$ of Coxeter groups.
Moreover, by Proposition~\ref{exactness of Enright} the functor
$S=S_{i_1}(a_1)\circ\cdots\circ S_{i_r}(a_r)$ with $a_k=\lan
h_{i_k},s_{i_{k+1}}\cdots s_{i_r}\circ\lam\ran$ induces a category
equivalence $\BM_i(\Gg,\lan h_i,\lam\ran)\to\BM_i(\Gg,\lan h_i,\lam'\ran)$
and we have $S(M(w\circ\lam))=M(xwx^{-1}\circ\lam')$,
$S(L(w\circ\lam))=L(xwx^{-1}\circ\lam')$.
Thus the proof of Theorem~\ref{Main theorem} in the case
$\BQ\Delta(\lam)\not\ni\delta$ is reduced to the case where
$\Delta(\lam)\subset\Delta_J$ for a proper subset $J$ of $I$.

Set
$$
\Gl_J=\Gh\oplus(\bigoplus_{\alpha\in\Delta_J}\Gg_\alpha),\quad
\Gn_J^+=\bigoplus_{\alpha\in\Delta^+\setminus\Delta_J}\Gg_\alpha,\quad
\Gn_J^-=\bigoplus_{\alpha\in\Delta^+\setminus\Delta_J}\Gg_{-\alpha},\quad
\Gp_J=\Gl_J\oplus\Gn_J^+.
$$
Note that we have $\dim\Gl_J<\infty$ since $J$ is a proper subset of $I$.
For $\mu\in\Gh^*$ let $M_J(\mu)$ be the Verma module for $\Gl_J$ with
highest weight $\mu$ and let $L_J(\mu)$ be its irreducible quotient.
We can regard them as $\Gp_J$-modules with trivial actions of $\Gn_J^+$.
By the definition we have $U(\Gg)\otimes_{U(\Gp_J)}M_J(\mu)\simeq M(\mu)$
for any $\mu\in\Gh^*$.
Hence Theorem~\ref{Main theorem} in the case $\BQ\Delta(\lam)\not\ni\delta$
follows from the character formula for the irreducible highest weight
modules over finite-dimensional semisimple Lie algebras, which is already
known (see the comments at the end), and the following result.
\begin{lemma}
For any $\lam\in\CC$ satisfying $\Delta(\lam)\subset\Delta_J$ we have
$U(\Gg)\otimes_{U(\Gp_J)}L_J(\lam)\simeq L(\lam)$.
\end{lemma}
\begin{proof}
Set $M=U(\Gg)\otimes_{U(\Gp_J)}L_J(\lam)$.
It is a highest weight module with highest weight $\lam$.
Set $M^{\Gn^+}=\{m\in M\,;\,\Gn^+m=0\}$.
It is sufficient to show $M^{\Gn^+}\cap M_{\lam-\xi}=0$ for any $\xi\in
Q^+\setminus\{0\}$.
Assume that $M^{\Gn^+}\cap M_{\lam-\xi}\ne\{0\}$ for some $\xi\in
Q^+\setminus\{0\}$.
By $\Delta(\lam)\subset\Delta_J$ and Proposition~\ref{prop:KK0} we have
$\xi\in\sum_{\alpha\in \Delta_J}\BZ\alpha$.
Hence under the isomorphism $M\simeq U(\Gn_J^-)\otimes_\BC L_J(\lam)$ we
have $M_{\lam-\xi}=1\otimes L_J(\lam)_{\lam-\xi}$.
It follows that $L_J(\lam)_{\lam-\xi}\cap L_J(\lam)^{\Gn^+\cap\Gl_J}\ne\{0\}$.
This contradicts the irreducibility of $L_J(\lam)$.
\end{proof}

The proof of Theorem~\ref{Main theorem} is complete in the case
$\BQ\Delta(\lam)\not\ni\delta$.

\bigskip

We finally give comments on the proof of the character formula for the
irreducible highest weight modules over finite-dimensional semisimple Lie
algebras which we have used in our proof in Case 2.
The unpublished result in the rational highest weight case due to
Beilinson-Bernstein (in particular, the part relating some twisted
$D$-modules with the twisted intersection cohomology groups of the Schubert
varieties) is recovered as a special case of the result in
Kashiwara-Tanisaki~\cite{KTneg2} (and also of the result in
Kashiwara-Tanisaki~\cite{KTpos3}).

\bibliographystyle{unsrt}
\def\same{\,$\raise3pt\hbox to 25pt{\hrulefill}\,$}

\end{document}